# SCALING LIMITS OF $(1+1)$-DIMENSIONAL PINNING MODELS WITH LAPLACIAN INTERACTION

BY FRANCESCO CARAVENNA[1] AND JEAN-DOMINIQUE DEUSCHEL

*Università degli Studi di Padova and TU Berlin*

We consider a random field $\varphi:\{1,\ldots,N\} \to \mathbb{R}$ with Laplacian interaction of the form $\sum_i V(\Delta\varphi_i)$, where $\Delta$ is the discrete Laplacian and the potential $V(\cdot)$ is symmetric and uniformly strictly convex. The *pinning model* is defined by giving the field a reward $\varepsilon \geq 0$ each time it touches the $x$-axis, that plays the role of a *defect line*. It is known that this model exhibits a phase transition between a delocalized regime ($\varepsilon < \varepsilon_c$) and a localized one ($\varepsilon > \varepsilon_c$), where $0 < \varepsilon_c < \infty$. In this paper we give a precise pathwise description of the transition, extracting the full scaling limits of the model. We show, in particular, that in the delocalized regime the field wanders away from the defect line at a typical distance $N^{3/2}$, while in the localized regime the distance is just $O((\log N)^2)$. A subtle scenario shows up in the critical regime ($\varepsilon = \varepsilon_c$), where the field, suitably rescaled, converges in distribution toward the derivative of a symmetric stable Lévy process of index 2/5. Our approach is based on Markov renewal theory.

## 1. Introduction.

1.1. *The model.* The main ingredient of our model is a function $V(\cdot):\mathbb{R} \to \mathbb{R} \cup \{+\infty\}$, that we call the *potential*. Our assumptions on $V(\cdot)$ are the following:

- *Symmetry*: $V(x) = V(-x)$, $\forall x \in \mathbb{R}$.
- *Uniform strict convexity*: there exists $\gamma > 0$ such that $x \mapsto V(x) - \gamma x^2/2$ is convex.
- *Regularity*: since $V(\cdot)$ is symmetric and convex, it is continuous and finite on some maximal interval $(-a, a)$ (possibly $a = +\infty$). We assume that

---

Received February 2008.

[1]Supported in part by the German Research Foundation–Research Group 718 during his stay at TU Berlin in January 2007.

$a > 0$ and we further require that $V(x) \to +\infty$ as $x \to \pm a$, so that the function $x \mapsto \exp(-V(x))$ is continuous on the whole real line.

Notice that, if $V(\cdot)$ is of class $C^2$ on $(-a, a)$, the uniform strict convexity assumption amounts to requiring that

$$\gamma := \inf_{x \in (-a,a)} V''(x) > 0. \tag{1.1}$$

It follows from the above assumptions that $\int_{\mathbb{R}} \exp(-V(x))\, dx < \infty$. Since adding a global constant to $V(\cdot)$ is immaterial for our purposes, we impose the normalization $\int_{\mathbb{R}} \exp(-V(x))\, dx = 1$. In this way we can interpret $\exp(-V(x))$ as a probability density, that has zero mean (by symmetry) and finite variance:

$$\sigma^2 := \int_{\mathbb{R}} x^2 e^{-V(x)}\, dx < \infty. \tag{1.2}$$

The most important example is, of course, the *Gaussian case*: $V(x) = x^2/2\sigma^2 + \log(\sigma\sqrt{2\pi})$.

Next we define the Hamiltonian, by setting

$$\mathcal{H}_{[a,b]}(\varphi) := \sum_{n=a+1}^{b-1} V(\Delta\varphi_n) \tag{1.3}$$

for $a, b \in \mathbb{Z}$ with $b - a \geq 2$ and for $\varphi: \{a, \ldots, b\} \to \mathbb{R}$, where $\Delta$ is the discrete Laplacian:

$$\Delta\varphi_n := (\varphi_{n+1} - \varphi_n) - (\varphi_n - \varphi_{n-1}) = \varphi_{n+1} + \varphi_{n-1} - 2\varphi_n. \tag{1.4}$$

We can now define our model: given $N \in \mathbb{N} := \{1, 2, \ldots\}$ and $\varepsilon \geq 0$, we introduce the probability measure $\mathbb{P}_{\varepsilon,N}$ on $\mathbb{R}^{N-1}$ defined by

$$\mathbb{P}_{\varepsilon,N}(d\varphi_1 \cdots d\varphi_{N-1}) = \frac{\exp(-\mathcal{H}_{[-1,N+1]}(\varphi))}{\mathcal{Z}_{\varepsilon,N}} \prod_{i=1}^{N-1} (\varepsilon \delta_0(d\varphi_i) + d\varphi_i), \tag{1.5}$$

where $d\varphi_i$ is the Lebesgue measure on $\mathbb{R}$, $\delta_0(\cdot)$ is the Dirac mass at zero and $\mathcal{Z}_{\varepsilon,N}$ is the normalization constant (*partition function*). To complete the definition, in order to make sense of $\mathcal{H}_{[-1,N+1]}(\varphi) = \mathcal{H}_{[-1,N+1]}(\varphi_{-1}, \varphi_0, \varphi_1, \ldots, \varphi_{N-1}, \varphi_N, \varphi_{N+1})$, we have to specify

(1.6)     the boundary conditions $\varphi_{-1} = \varphi_0 = \varphi_N = \varphi_{N+1} = 0$.

The choice of zero boundary conditions is made only for simplicity, but our approach and results go through for general choices (provided they are, say, bounded in $N$).

1.2. *The phase transition.* The law $\mathbb{P}_{\varepsilon,N}$ is what is called a *pinning model* and can be viewed as a $(1+1)$-dimensional model for a linear chain of length $N$ attracted to a defect line, namely, the $x$-axis. The parameter $\varepsilon \geq 0$ tunes the strength of the attraction and one wishes to understand its effect on the field, in the large $N$ limit.

The basic properties of this model (and of the closely related *wetting model*, in which the field is also constrained to stay nonnegative) were investigated in a first paper [6], to which we refer for a detailed discussion and for a survey of the literature. In particular, it was shown that there is a critical threshold $0 < \varepsilon_c < \infty$ that determines a *phase transition* between a *delocalized regime* ($\varepsilon < \varepsilon_c$), in which the reward is essentially ineffective, and a *localized regime* ($\varepsilon > \varepsilon_c$), in which on the other hand the reward has a macroscopic effect on the field. More precisely, defining the *contact number* $\ell_N$ by

$$(1.7) \qquad \ell_N := \#\{i \in \{1,\ldots,N\} : \varphi_i = 0\},$$

we have the following dichotomy:

- if $\varepsilon \leq \varepsilon_c$, then for every $\delta > 0$ and $N \in \mathbb{N}$

$$(1.8) \qquad \mathbb{P}_{\varepsilon,N}\left(\frac{\ell_N}{N} > \delta\right) \leq e^{-c_1 N},$$

  where $c_1$ is a positive constant;
- if $\varepsilon > \varepsilon_c$, then there exists $\mathrm{D}(\varepsilon) > 0$ such that for every $\delta > 0$ and $N \in \mathbb{N}$

$$(1.9) \qquad \mathbb{P}_{\varepsilon,N}\left(\left|\frac{\ell_N}{N} - \mathrm{D}(\varepsilon)\right| > \delta\right) \leq e^{-c_2 N},$$

  where $c_2$ is a positive constant.

Roughly speaking, for $\varepsilon \leq \varepsilon_c$ we have $\ell_N = o(N)$, while for $\varepsilon > \varepsilon_c$ we have $\ell_N \sim \mathrm{D}(\varepsilon) \cdot N$. For an explicit characterization of $\varepsilon_c$ and $\mathrm{D}(\varepsilon)$ we refer to [6], where it is also proven that the phase transition is exactly of second order.

The aim of this paper is to go far beyond (1.8) and (1.9) in the study of the path properties of $\mathbb{P}_{\varepsilon,N}$. Our results, that include the scaling limits of the model on $C([0,1])$, provide strong path characterizations of (de)localization. We also show that the delocalized regime ($\varepsilon < \varepsilon_c$) and the critical one ($\varepsilon = \varepsilon_c$) exhibit great differences, that are somewhat hidden in relation (1.8). In fact, a closer look at the critical regime exposes a rich structure that we analyze in detail.

REMARK 1.1. We point out that the hypothesis on $V(\cdot)$ in the present paper are stronger than those of [6] [where essentially only the second moment condition (1.2) was required]. This is a price to pay in order to obtain precise path results, like, for instance, Theorem 1.4 below, that would not hold in the general setting of [6]. Although for some other results our assumptions could have been weakened, we have decided not to do it, both to keep us in a unified setting, and because, with the uniform strict convexity assumption on $V(\cdot)$, one can apply general powerful tools, notably the Brascamp–Lieb inequality [4, 5], that allow to give streamlined versions of otherwise rather technical proofs.

Also notice that the analysis of this paper does not cover the wetting model, that was also considered in [6]. The reason for this exclusion is twofold: on the one hand, the basic estimates derived in [6] in the wetting case are not sufficiently precise as those obtained in the pinning case; on the other hand, for the scaling limits of the wetting model one should rely on suitable invariance principles for the integrated random walk process conditioned to stay nonnegative, but this issue seems not to have been investigated in the literature.

1.3. *Path results and the scaling limits.* Let us look first at the free case $\varepsilon = 0$, where the pinning reward is absent. It was shown in [6] that the law $\mathbb{P}_{0,N}$ enjoys the following random walk interpretation (for more details see Section 2). Let $(\{Y_n\}_{n \in \mathbb{Z}^+ := \mathbb{N} \cup \{0\}}, \mathbf{P})$ denote a real random walk starting at zero and with step law $\mathbf{P}(Y_1 \in dx) = \exp(-V(x)) dx$, and let $\{Z_n\}_{n \in \mathbb{Z}^+}$ be the corresponding *integrated random walk process*:

$$Z_0 := 0, \qquad Z_n := Y_1 + \cdots + Y_n, \qquad n \in \mathbb{N}.$$

The basic fact is that $\mathbb{P}_{0,N}$ coincides with the law of the vector $(Z_1, \ldots, Z_{N-1})$ conditionally on $(Z_N, Z_{N+1}) = (0,0)$, that is, the free law $\mathbb{P}_{0,N}$ is nothing but the bridge of an integrated random walk. By our assumptions on $V(\cdot)$ (see Section 1.1), the walk $\{Y_n\}_n$ has zero mean and finite variance $\sigma^2$, hence, for large $k$ the variable $Z_k$ scales like $k^{3/2}$. It is therefore natural to consider the following rescaled and linearly-interpolated version of the field $\{\varphi_n\}_n$:

$$\widehat{\varphi}_N(t) := \frac{\varphi_{\lfloor Nt \rfloor}}{\sigma N^{3/2}} + (Nt - \lfloor Nt \rfloor) \frac{\varphi_{\lfloor Nt \rfloor + 1} - \varphi_{\lfloor Nt \rfloor}}{\sigma N^{3/2}},$$

(1.10)
$$N \in \mathbb{N}, t \in [0,1],$$

and to study the convergence in distribution of $\{\widehat{\varphi}_N(t)\}_{t \in [0,1]}$ as $N \to \infty$ on $C([0,1])$, the space of real valued, continuous functions on $[0,1]$ (equipped as usual with the topology of uniform convergence). To this purpose, we let $\{B_t\}_{t \in [0,1]}$ denote a standard Brownian motion on the interval $[0,1]$, we define the *integrated Brownian motion process* $\{I_t\}_{t \in [0,1]}$ by $I_t := \int_0^t B_s \, ds$ and we introduce the conditioned process

(1.11) $\{(\widehat{B}_t, \widehat{I}_t)\}_{t \in [0,1]} := \{(B_t, I_t)\}_{t \in [0,1]}$ conditionally on $(B_1, I_1) = (0,0)$.

Exploiting the random walk description of $\mathbb{P}_{0,N}$, it is not difficult to show that the process $\{\widehat{\varphi}_N(t)\}_{t \in [0,1]}$ under $\mathbb{P}_{0,N}$ converges in distribution as $N \to \infty$ toward $\{\widehat{I}_t\}_{t \in [0,1]}$. The emergence of a nontrivial scaling limit for $\{\widehat{\varphi}_N(t)\}_{t \in [0,1]}$ is a precise formulation of the statement that the typical height of the field under $\mathbb{P}_{0,N}$ is of order $N^{3/2}$. It is natural to wonder what happens of this picture when $\varepsilon > 0$: the answer is given by our first result.

THEOREM 1.2 (Scaling limits). *The rescaled field $\{\widehat{\varphi}_N(t)\}_{t \in [0,1]}$ under $\mathbb{P}_{\varepsilon,N}$ converges in distribution on $C([0,1])$ as $N \to \infty$, for every $\varepsilon \geq 0$. The limit is as follows:*

- *If $\varepsilon < \varepsilon_c$, the law of the process $\{\widehat{I}_t\}_{t \in [0,1]}$;*
- *If $\varepsilon = \varepsilon_c$ or $\varepsilon > \varepsilon_c$, the law concentrated on the constant function $f(t) \equiv 0$, $t \in [0,1]$.*

Thus, the pinning reward $\varepsilon$ is ineffective for $\varepsilon < \varepsilon_c$, at least for the large scale properties of the field that are identical to the free case $\varepsilon = 0$. On the other hand, if $\varepsilon \geq \varepsilon_c$, the reward is able to change the macroscopic behavior of the field, whose height under $\mathbb{P}_{\varepsilon,N}$ scales less than $N^{3/2}$. We are now going to strengthen these considerations by looking at path properties on a finer scale. However, before proceeding, we stress that, from the point of view of the scaling limits, the critical regime $\varepsilon = \varepsilon_c$ is close to the localized one $\varepsilon > \varepsilon_c$ rather than to the delocalized one $\varepsilon < \varepsilon_c$, in contrast with (1.8) and (1.9).

We start looking at the delocalized regime ($\varepsilon < \varepsilon_c$). It is convenient to introduce the *contact set* $\tau$ of the field $\{\varphi_i\}_{i \in \mathbb{Z}^+}$, that is, the random subset of $\mathbb{Z}^+$ defined by

$$(1.12) \qquad \tau := \{i \in \mathbb{Z}^+ : \varphi_i = 0\} \subseteq \mathbb{Z}^+,$$

where we set by definition $\varphi_0 := 0$ so that $0 \in \tau$. We already know from (1.8) that for $\varepsilon < \varepsilon_c$ we have $\#\{\tau \cap [0,N]\} = \ell_N + 1 = o(N)$ under $\mathbb{P}_{\varepsilon,N}$. The next theorem shows that in fact $\tau \cap [0,N]$ consists of a finite number of points (i.e., the variable $\ell_N$ under $\mathbb{P}_{\varepsilon,N}$ is tight) and all these points are at finite distance from the boundary.

THEOREM 1.3. *For every $\varepsilon < \varepsilon_c$ the following relation holds:*

$$(1.13) \qquad \lim_{L \to \infty} \liminf_{N \to \infty} \mathbb{P}_{\varepsilon,N}(\tau \cap [L, N-L] = \varnothing) = 1.$$

We will see that the scaling limit of $\{\widehat{\varphi}_N(t)\}_{t \in [0,1]}$ under $\mathbb{P}_{\varepsilon,N}$, for $\varepsilon \in (0, \varepsilon_c)$, is a direct consequence of relation (1.13) and of the scaling limit for $\varepsilon = 0$. The reason for this lies in the following crucial fact: conditionally on the contact set, the excursion of the field under $\mathbb{P}_{\varepsilon,N}$ between two consecutive contact points, say, $\tau_k$ and $\tau_{k+1}$, is distributed according to the free law $\mathbb{P}_{0,\tau_{k+1}-\tau_k}$ with suitable boundary conditions (see Section 2.3 for more details).

Next we focus on the localized and critical regimes $(\varepsilon > \varepsilon_c)$ and $(\varepsilon = \varepsilon_c)$. The first question left open by Theorem 1.2 is, of course, if one can obtain a more precise estimate on the height of the field than just $o(N^{3/2})$. We have the following result.

THEOREM 1.4. *For every $\varepsilon > \varepsilon_c$ the following relation holds:*

$$\lim_{K \to \infty} \liminf_{N \to \infty} \mathbb{P}_{\varepsilon,N}\left(\max_{0 \leq k \leq N} |\varphi_k| \leq K(\log N)^2\right) = 1 \tag{1.14}$$

*while for $\varepsilon = \varepsilon_c$ the following relation holds:*

$$\lim_{K \to \infty} \liminf_{N \to \infty} \mathbb{P}_{\varepsilon_c,N}\left(\frac{1}{K}\frac{N^{3/2}}{(\log N)^{3/2}} \leq \max_{0 \leq k \leq N} |\varphi_k| \leq K\frac{N^{3/2}}{\log N}\right) = 1. \tag{1.15}$$

The fact that $\{\widehat{\varphi}_N(t)\}_{t \in [0,1]}$ under $\mathbb{P}_{\varepsilon,N}$ has, for $\varepsilon \geq \varepsilon_c$, a trivial scaling limit, is of course an immediate consequence of the upper bounds on $\max_{0 \leq k \leq N} |\varphi_k|$ in (1.14) and (1.15).

We believe that the optimal scaling of $\max_{0 \leq k \leq N} |\varphi_k|$ for $\varepsilon = \varepsilon_c$ is given by the lower bound in (1.15) (to lighten the exposition, we do not investigate this problem deeper).

REMARK 1.5. Another interesting quantity is the maximal gap $\Delta_N$, defined as

$$\Delta_N := \max\{\tau_k - \tau_{k-1} : 0 \leq k \leq \ell_N\}. \tag{1.16}$$

We already know from (1.13) that $\Delta_N \sim N$ in the delocalized regime $(\varepsilon < \varepsilon_c)$. It turns out that in the localized regime $(\varepsilon > \varepsilon_c)$ we have $\Delta_N = O(\log N)$, while in the critical regime $(\varepsilon = \varepsilon_c)$ $\Delta_N \approx N/\log N$, meaning by this that

$$\lim_{K \to \infty} \liminf_{N \to \infty} \mathbb{P}_{\varepsilon_c,N}\left(\frac{1}{K}\frac{N}{\log N} \leq \Delta_N \leq K\frac{N}{\log N}\right) = 1.$$

For $\varepsilon = \varepsilon_c$ we also have $\ell_N \approx N/\log N$. For conciseness, we omit a detailed proof of these relations (though some partial results will be given in the proof of Theorem 1.4, see also Appendix A). Table 1 summarizes the results described so far.

1.4. *A refined critical scaling limit.* Relation (1.15) shows that the field in the critical regime has very large fluctuations, *almost* of the order $N^{3/2}$. This may suggest the possibility of lowering the scaling constants $N^{3/2}$ in the definition (1.10) of the rescaled field $\widehat{\varphi}_N(t)$, in order to make a nontrivial scaling limit emerge under $\mathbb{P}_{\varepsilon_c,N}$. However, some care is needed: in fact, $\Delta_N/N \to 0$ as $N \to \infty$ under $\mathbb{P}_{\varepsilon_c,N}$ and this means that, *independently of the choice of the scaling constants*, the zero level set of the rescaled field becomes dense in $[0,1]$. This fact rules out the possibility of getting a nontrivial scaling limit in $C([0,1])$, or even in the space of càdlàg functions $D([0,1])$.

TABLE 1
*A schematic representation of the order of growth as $N \to \infty$ of the three quantities $\max_{0 \le k \le N} |\varphi_k|$, $\Delta_N$ and $\ell_N$ under $\mathbb{P}_{\varepsilon,N}$*

|  | $(\varepsilon < \varepsilon_c)$ | $(\varepsilon = \varepsilon_c)$ | $(\varepsilon > \varepsilon_c)$ |
|---|---|---|---|
| $\max_{0 \le k \le N} |\varphi_k|$ | $N^{3/2}$ | $\frac{N^{3/2}}{(\log N)^{3/2}} \div \frac{N^{3/2}}{\log N}$ | $O((\log N)^2)$ |
| $\Delta_N$ | $N$ | $\frac{N}{\log N}$ | $O(\log N)$ |
| $\ell_N$ | $O(1)$ | $\frac{N}{\log N}$ | $N$ |

We are going to show that a nontrivial scaling limit can indeed be extracted in a distributional sense, that is, integrating the field against test functions, and to this purpose, the right scaling constants turn out to be $N^{3/2}/(\log N)^{5/2}$ (see below for an heuristic explanation). Therefore, we introduce the new rescaled field $\{\widetilde{\varphi}_N(t)\}_{t \in [0,1]}$ (this time with no need of linear interpolation) defined by

$$(1.17) \qquad \widetilde{\varphi}_N(t) := \frac{(\log N)^{5/2}}{N^{3/2}} \varphi_{\lfloor Nt \rfloor}.$$

Viewing $\widetilde{\varphi}_N(t)$ as a density, we introduce the signed measure $\boldsymbol{\mu}_N$ on $[0,1]$ defined by

$$(1.18) \qquad \boldsymbol{\mu}_N(\mathrm{d}t) := \widetilde{\varphi}_N(t) \, \mathrm{d}t.$$

We look at $\boldsymbol{\mu}_N$ under the critical law $\mathbb{P}_{\varepsilon_c,N}$ as a random element of $\mathcal{M}([0,1])$, the space of all *finite signed Borel measures* on the interval $[0,1]$, that we equip with the topology of vague convergence and with the corresponding Borel $\sigma$-field ($\nu_n \to \nu$ vaguely if and only if $\int f \, \mathrm{d}\nu_n \to \int f \, \mathrm{d}\nu$ for all bounded and continuous functions $f:[0,1] \to \mathbb{R}$). Our goal is to show that the sequence $\{\boldsymbol{\mu}_N\}_N$ has a nontrivial limit in distribution on $\mathcal{M}([0,1])$.

To describe the limit, let $\{L_t\}_{t \ge 0}$ denote the stable symmetric Lévy process of index $2/5$ (a standard version with càdlàg paths). More explicitly, $\{L_t\}_{t \ge 0}$ is a Lévy process with zero drift, zero Brownian component and with Lévy measure given by $\Pi(\mathrm{d}x) = c_L |x|^{-7/5} \, \mathrm{d}x$, where the positive constant $c_L$ is defined explicitly in equation (6.25). Since the index is less than 1, the paths of $L$ are a.s. of bounded variation (cf. [2]), hence, we can define path by path the (random) finite signed measure $\mathrm{d}L$ in the Steltjes sense:

$$\mathrm{d}L((a,b]) := L_b - L_a.$$

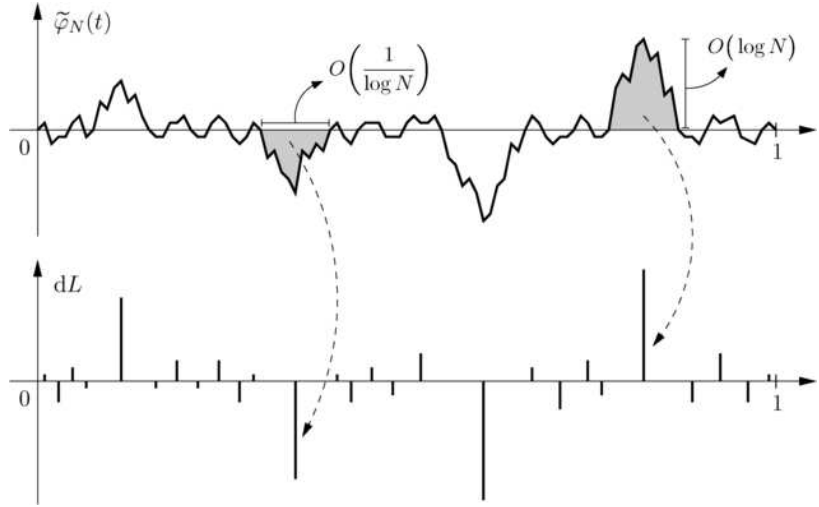

FIG. 1. *A graphical representation of Theorem 1.6. For large $N$, the excursions of the rescaled field under the critical law $\mathbb{P}_{\varepsilon_c,N}$ contribute to the measure $\boldsymbol{\mu}_N(\mathrm{d}t)$ [see (1.18)], approximately like Dirac masses, with intensity given by their (signed) area. The width and height of the relevant excursions are of order $(1/\log N)$ and $\log N$ respectively. We warn the reader that the $x$- and $y$-axis in the picture have different units of length, and that the field can actually cross the $x$-axis without touching it (though this feature has not been evidenced in the picture for simplicity).*

We stress that $\mathrm{d}L$ is a.s. a purely atomic measure, that is, a sum of Dirac masses (for more details and for an explicit construction of $\mathrm{d}L$, see Remark 1.7 below).

We are now ready to state our main result (see Figure 1 for a graphical description).

THEOREM 1.6. *The random signed measure $\boldsymbol{\mu}_N$ under $\mathbb{P}_{\varepsilon_c,N}$ converges in distribution on $\mathcal{M}([0,1])$ as $N \to \infty$ toward the the random signed measure $\mathrm{d}L$.*

This result describes in a quantitative way the rich structure of the field for $\varepsilon = \varepsilon_c$. Let us try to give a heuristic description. Roughly speaking, for large $N$ the profile of the unrescaled field $\{\varphi_i\}_{0 \le i \le N}$ under $\mathbb{P}_{\varepsilon_c,N}$ is dominated by the *large excursions* over the contact set, that is, by those excursions whose width is of the same order $\approx N/\log N$ as the maximal gap $\Delta_N$ (see Table 1). As already observed, each excursion is distributed according to the free law (with suitable boundary conditions), hence, by Theorem 1.2, the height of these excursions is of order $\approx (N/\log N)^{3/2}$. When the field is rescaled according to (1.17), the width of these excursions becomes of order $\approx 1/\log N$ and their height of order $\approx \log N$, hence, for large $N$ they contribute to the measure $\boldsymbol{\mu}_N$ approximately like Dirac masses (see Figure 1). Therefore, the properties of these large excursions, for large $N$, can be read from the structure of the Dirac masses that build the limit measure $dL$; see (1.19) below.

REMARK 1.7. The measure $dL$ can be constructed in the following explicit way; compare [2]. Let $\mathcal{S}$ denote a *Poisson point process* on the space $\mathsf{X} := [0,1] \times \mathbb{R}$ with intensity measure $\gamma := dx \otimes c_L |y|^{-7/5} dy$ (where $dx$ and $dy$ denote the Lebesgue measure on $[0,1]$ and $\mathbb{R}$). We recall that $\mathcal{S}$ is a random countable subset of $\mathsf{X}$ with the following properties:

— for every Borel set $A \subseteq \mathsf{X}$, the random variable $\#(\mathcal{S} \cap A)$ has a Poisson distribution with parameter $\gamma(A)$ [the symbol $\#$ denotes the cardinality of a set and in case $\gamma(A) = +\infty$ we mean that $\#(\mathcal{S} \cap A) = +\infty$, a.s.];
— for any $k \in \mathbb{N}$ and for every family of pairwise disjoint Borel sets $A_1, \ldots, A_k \subseteq \mathsf{X}$, the random variables $\#(\mathcal{S} \cap A_1), \ldots, \#(\mathcal{S} \cap A_k)$ are independent.

Since $\gamma$ is a $\sigma$-finite measure, the random set $\mathcal{S}$ is a.s. countable: enumerating its points in some (arbitrary) way, say, $\mathcal{S} = \{(x_i, y_i)\}_{i \in \mathbb{N}}$, we can write

$$(1.19) \qquad dL(\cdot) \stackrel{d}{=} \sum_{i \in \mathbb{N}} y_i \cdot \delta_{x_i}(\cdot),$$

where $\delta_x(\cdot)$ denotes the Dirac mass at $x \in \mathbb{R}$. Notice that since $\int_{\mathsf{X}} (|y| \wedge 1) d\gamma < \infty$, the r.h.s. of (1.19) is indeed a finite measure, that is, $\sum_{i \in \mathbb{N}} |y_i| < \infty$ a.s.; compare [11].

1.5. *Outline of the paper.* The exposition is organized as follows:

- In Section 2 we recall some basic properties of $\mathbb{P}_{\varepsilon,N}$ that have been proven in [6]. In particular, we develop a renewal theory description of the model, which is the cornerstone of our approach.
- In Section 3 we prove a part of Theorem 1.4, more precisely, equation (1.14) and the upper bound on $\max_{0 \leq k \leq N} |\varphi_k|$ in (1.15), exploiting the Brascamp–Lieb inequality. These results also prove Theorem 1.2 for $\varepsilon \geq \varepsilon_c$.
- Section 4 is devoted to the proof of Theorem 1.3 and of Theorem 1.2 for $\varepsilon < \varepsilon_c$.
- In Section 5 we complete the proof of Theorem 1.4, obtaining the lower bound on $\max_{0 \leq k \leq N} |\varphi_k|$ in equation (1.15).
- Section 6 is devoted to the proof of Theorem 1.6.
- Finally, some technical points are treated in the Appendixes A and B.

**2. Some basic facts.** This section is devoted to a detailed description of $\mathbb{P}_{\varepsilon,N}$, taking inspiration from [6]. We show in Section 2.1 that, conditionally on the contact set $\tau$ [cf. (1.12)], the pinning model $\mathbb{P}_{\varepsilon,N}$ is linked to the integral of a random walk. Then in Section 2.2 we focus on the law of the contact set itself, which admits a crucial description in terms of Markov renewal theory. We conclude by putting together these results in Section 2.3, where we show that the full measure $\mathbb{P}_{\varepsilon,N}$ is the conditioning of an explicit infinite-volume law $\mathcal{P}_\varepsilon$.

2.1. *Integrated random walk.* One of the key features of the model $\mathbb{P}_{\varepsilon,N}$ is its link with the integral of a random walk, described in Section 2 of [6], that we now recall.

Given $a, b \in \mathbb{R}$, we define on some probability space $(\Omega, \mathcal{F}, \mathbf{P} = \mathbf{P}^{(a,b)})$ a sequence $\{X_i\}_{i \in \mathbb{N}}$ of independent and identically distributed random variables, with marginal laws $X_1 \sim \exp(-V(x))\,\mathrm{d}x$. By our assumptions on $V(\cdot)$, it follows that

$$\mathbf{E}(X_1) = 0, \qquad \mathbf{E}(X_1^2) = \sigma^2 \in (0, \infty).$$

We denote by $\{Y_i\}_{i \in \mathbb{Z}^+}$ the associated random walk starting at $a$, that is,

(2.1) $$Y_0 = a, \qquad Y_n = a + X_1 + \cdots + X_n, \qquad n \in \mathbb{N},$$

while $\{Z_i\}_{i \in \mathbb{Z}^+}$ denotes the *integrated random walk* starting at $b$, that is, $Z_0 = b$ and for $n \in \mathbb{N}$

(2.2)
$$\begin{aligned} Z_n &= b + Y_1 + \cdots + Y_n \\ &= b + na + nX_1 + (n-1)X_2 + \cdots + 2X_{n-1} + X_n. \end{aligned}$$

Notice that

(2.3) $$\{(Y_n, Z_n)\}_n \quad \text{under} \quad \mathbf{P}^{(a,b)} \stackrel{d}{=} \{(Y_n + a, Z_n + b + na)\}_n \quad \text{under} \quad \mathbf{P}^{(0,0)}.$$

The marginal distributions of the process $\{Z_n\}_n$ are easily computed [6], Lemma 4.2:

(2.4)
$$\begin{aligned} &\mathbf{P}^{(a,b)}((Z_1, \ldots, Z_n) \in (\mathrm{d}z_1, \ldots, \mathrm{d}z_n)) \\ &\qquad = e^{-\mathcal{H}_{[-1,n]}(b-a,b,z_1,\ldots,z_n)}\,\mathrm{d}z_1 \cdots \mathrm{d}z_n, \end{aligned}$$

where $\mathcal{H}_{[-1,n]}(\cdot)$ is exactly the Hamiltonian of our model, defined in (1.3).

We are ready to make the link with our model $\mathbb{P}_{\varepsilon,N}$. In the free case $\varepsilon = 0$, it is rather clear from (2.4) and (1.5) that $\mathbb{P}_{0,N}$ is nothing but the law of $(Z_1, \ldots, Z_{N-1})$ under $\mathbf{P}^{(0,0)}(\cdot | Z_N = 0, Z_{N+1} = 0)$, that is, the polymer measure $\mathbb{P}_{0,N}$ is just the bridge of the integral of a random walk; compare [6], Proposition 4.1.

To deal with the case $\varepsilon > 0$, we recall the definition of the contact set $\tau$, given in (1.12):

$$\tau := \{i \in \mathbb{Z}^+ : \varphi = 0\} \subseteq \mathbb{Z}^+.$$

We also set for conciseness $\tau_{[a,b]} := \tau \cap [a,b]$. Then, again comparing (2.4) with (1.5), we have the following basic relation: for $\varepsilon > 0$, $N \in \mathbb{N}$ and for every subset $A \subseteq \{1, \ldots, N-1\}$,

$$
(2.5) \quad \begin{aligned} \mathbb{P}_{\varepsilon,N}(\cdot | \tau_{[1,N-1]} = A) \\ = \mathbf{P}^{(0,0)}((Z_1, \ldots, Z_{N-1}) \in \cdot | Z_i = 0, \ \forall i \in A \cup \{N, N+1\}). \end{aligned}
$$

In words, once we fix the contact set $\tau_{[1,N-1]} = A$, the field $(\varphi_1, \ldots, \varphi_{N-1})$ under $\mathbb{P}_{\varepsilon,N}$ is distributed like the integrated random walk $(Z_1, \ldots, Z_{N-1})$ under $\mathbf{P}^{(0,0)}$ conditioned on being zero at the epochs in $A$ and also at $N$ and $N+1$ [because of the boundary conditions, cf. (1.6)]. A crucial aspect of (2.5) is that the r.h.s. is independent of $\varepsilon$. Therefore, all the dependence of $\varepsilon$ of $\mathbb{P}_{\varepsilon,N}$ is contained in the law of the contact set.

Notice that in the l.h.s. of (2.5) we are really conditioning on an event of positive probability, while the conditioning in the r.h.s. of (2.5) is to be understood in the sense of conditional distributions [which can be defined unambiguously, because we have assumed that the density $x \mapsto e^{-V(x)}$ is continuous].

We conclude this paragraph observing that the joint process $\{(Y_n, Z_n)\}_n$ under $\mathbf{P}^{(a,b)}$ is a Markov process on $\mathbb{R}^2$. On the other hand, the process $\{Z_n\}_n$ alone is not Markov, but it rather has *finite memory* $m = 2$, that is, for every $n \in \mathbb{N}$

$$
(2.6) \quad \begin{aligned} \mathbf{P}^{(a,b)}(\{Z_{n+k}\}_{k \geq 1} \in \cdot | Z_i, i \leq n) &= \mathbf{P}^{(a,b)}(\{Z_{n+k}\}_{k \geq 1} \in \cdot | Z_{n-1}, Z_n) \\ &= \mathbf{P}^{(Z_n - Z_{n-1}, Z_n)}(\{Z_k\}_{k \geq 1} \in \cdot). \end{aligned}
$$

In fact, from (2.4) it is clear that

$$(2.7) \qquad \mathbf{P}^{(a,b)} = \mathbf{P}(\cdot | Z_{-1} = b - a, Z_0 = b).$$

2.2. *Markov renewal theory.* It is convenient to identify the random set $\tau$ with the increasing sequence of variables $\{\tau_k\}_{k \in \mathbb{Z}^+}$ defined by

$$(2.8) \qquad \tau_0 := 0, \qquad \tau_{k+1} := \inf\{i > \tau_k : \varphi_i = 0\}.$$

Observe that the contact number $\ell_N$, introduced in (1.7), can be also expressed as

$$(2.9) \qquad \ell_N = \max\{k \in \mathbb{Z}^+ : \tau_k \leq N\} = \sum_{k=1}^{N} \mathbf{1}_{\{k \in \tau\}}.$$

We also introduce the process $\{J_k\}_{k \in \mathbb{Z}^+}$ that gives the height of the field just before the contact points:

$$(2.10) \qquad J_0 := 0, \qquad J_k := \varphi_{\tau_k - 1}, \qquad k \in \mathbb{N}.$$

Of course, under the law $\mathbb{P}_{\varepsilon,N}$, we look at the variables $\tau_k$, $J_k$ only for $k \leq \ell_N$. The crucial fact, proven in Section 3 of [6], is that the vector $\{\ell_N, (\tau_k)_{k \leq \ell_N}, (J_k)_{k \leq \ell_N}\}$ under the law $\mathbb{P}_{\varepsilon,N}$ admits an explicit description in terms of Markov renewal theory, that we now recall.

Following [6], Section 3.2, for $\varepsilon > 0$ we denote by $\mathcal{P}_\varepsilon$ the law under which the joint process $\{(\tau_k, J_k)\}_{k \in \mathbb{Z}^+}$ is a Markov process on $\mathbb{Z}^+ \cup \{\infty\} \times \mathbb{R} \cup \{\infty\}$, with starting point $(\tau_0, J_0) = (0, 0)$ and with transition kernel given by

$$(2.11) \quad \mathcal{P}_\varepsilon((\tau_{k+1}, J_{k+1}) \in (\{n\}, \mathrm{d}y) | (\tau_k, J_k) = (m, x)) := \mathsf{K}^\varepsilon_{x,\mathrm{d}y}(n-m),$$

where $\mathsf{K}^\varepsilon_{x,\mathrm{d}y}(n)$ is defined for $x, y \in \mathbb{R}$ and $n \in \mathbb{N}$ by

$$(2.12) \quad \mathsf{K}^\varepsilon_{x,\mathrm{d}y}(n) := \varepsilon e^{-\mathrm{F}(\varepsilon)n} \frac{v_\varepsilon(y)}{v_\varepsilon(x)} \cdot \left. \frac{\mathbf{P}^{(-x,0)}(Z_{n-1} \in \mathrm{d}y, Z_n \in \mathrm{d}z)}{\mathrm{d}z} \right|_{z=0}.$$

The function $\mathrm{F}(\varepsilon)$ is the *free energy* of the model $\mathbb{P}_{\varepsilon,N}$, while $v_\varepsilon(\cdot)$ is a suitable positive function connected to an infinite dimensional Perron–Frobenius eigenvalue problem (we refer to Sections 3 and 4 of [6] for a detailed discussion). We stress that $\mathrm{F}(\varepsilon) = 0$ if $\varepsilon \leq \varepsilon_c$, while $\mathrm{F}(\varepsilon) > 0$ if $\varepsilon > \varepsilon_c$.

The dependence of the r.h.s. of (2.11) on $n - m$ implies that, under $\mathcal{P}_\varepsilon$, the process $\{J_k\}_k$ alone is itself a Markov chain on $\mathbb{R} \cup \{\infty\}$, with transition kernel

$$(2.13) \qquad \mathcal{P}_\varepsilon(J_{k+1} \in \mathrm{d}y | J_k = x) = D^\varepsilon_{x,\mathrm{d}y} := \sum_{n \in \mathbb{N}} \mathsf{K}^\varepsilon_{x,\mathrm{d}y}(n).$$

On the other hand, the process $\{\tau_k\}_k$ is not a Markov chain, but rather a *Markov renewal process* (cf. [1]) in fact, its increments $\{\tau_{k+1} - \tau_k\}_k$ are independent conditionally on the *modulating chain* $\{J_k\}_{k \in \mathbb{Z}^+}$, as it is clear from (2.11).

From (2.12) it follows that, as a measure in $\mathrm{d}y$, the kernel $\mathsf{K}^\varepsilon_{x,\mathrm{d}y}(n)$ is absolutely continuous for $n \geq 2$, while $\mathsf{K}^\varepsilon_{x,\mathrm{d}y}(1)$ is a multiple of the Dirac mass at zero $\delta_0(\mathrm{d}y)$. The properties of the kernel $\mathsf{K}^\varepsilon_{x,\mathrm{d}y}(n)$ depend strongly on the value of $\varepsilon$. First, the kernel is defective if $\varepsilon < \varepsilon_c$, while it is proper if $\varepsilon \geq \varepsilon_c$, since

$$\sum_{n \in \mathbb{N}} \int_{y \in \mathbb{R}} \mathsf{K}^\varepsilon_{x,\mathrm{d}y}(n) = \int_{y \in \mathbb{R}} D^\varepsilon_{x,\mathrm{d}y} = \min\left\{\frac{\varepsilon}{\varepsilon_c}, 1\right\},$$

so that the probability that $\tau_k = \infty$ for some $k$ is one if $\varepsilon < \varepsilon_c$, while it is zero if $\varepsilon \geq \varepsilon_c$. Moreover, as $n \to \infty$ for fixed $x, y \in \mathbb{R}$, we have

$$(2.14) \qquad \frac{\mathsf{K}^\varepsilon_{x,\mathrm{d}y}(n)}{\mathrm{d}y} = \frac{L^\varepsilon(x,y)}{n^2} e^{-\mathrm{F}(\varepsilon) \cdot n}(1 + o(1))$$

for a suitable function $L^\varepsilon(x, y)$, compare [6], Sections 3.2 and 4.1.

To summarize, the kernel $\mathsf{K}^\varepsilon_{x,\mathrm{d}y}(n)$ is defective with heavy tails in the delocalized regime ($\varepsilon < \varepsilon_c$) and it is proper with heavy tails in the critical regime ($\varepsilon = \varepsilon_c$), while it is proper with exponential tails in the localized regime ($\varepsilon > \varepsilon_c$). We also note that when $\varepsilon \geq \varepsilon_c$ the modulating chain $\{J_k\}_{k\in\mathbb{Z}^+}$ on $\mathbb{R}$ is *positive recurrent*, that is, it admits an invariant probability law $\nu_\varepsilon$: $\int_{x\in\mathbb{R}} \nu_\varepsilon(\mathrm{d}x) D^\varepsilon_{x,\mathrm{d}y} = \nu_\varepsilon(\mathrm{d}y)$, with $\nu_\varepsilon(\{0\}) > 0$ and no other atoms.

We are now ready to link the law $\mathcal{P}_\varepsilon$ to our model $\mathbb{P}_{\varepsilon,N}$. Introducing the event

$$\mathcal{A}_N := \{\{N, N+1\} \subseteq \tau\} = \{\tau_k = N, \tau_{k+1} = N+1 \text{ for some } k \in \mathbb{N}\},$$

Proposition 3.1 of [6] states that the vector $\{\ell_N, (\tau_k)_{k\leq\ell_N}, (J_k)_{k\leq\ell_N}\}$ has the same distribution under the law $\mathbb{P}_{\varepsilon,N}$ and under $\mathcal{P}_\varepsilon(\cdot|\mathcal{A}_N)$. More precisely, for all $k \in \mathbb{N}$, $\{t_i\}_{1\leq i\leq k} \in \mathbb{N}^k$ and $\{y_i\}_{1\leq i\leq k} \in \mathbb{R}^k$ we have

$$\begin{aligned}
\mathbb{P}_{\varepsilon,N}(\ell_N = k, \tau_i = t_i, J_i \in \mathrm{d}y_i, i \leq k) \\
= \mathcal{P}_\varepsilon(\ell_N = k, \tau_i = t_i, J_i \in \mathrm{d}y_i, i \leq k | \mathcal{A}_N).
\end{aligned} \quad (2.15)$$

In words, the contact set $\tau \cap [0, N]$ under the law $\mathbb{P}_{\varepsilon,N}$ is distributed like a Markov renewal process, of law $\mathcal{P}_\varepsilon$ and modulating chain $\{J_k\}_k$, conditioned to visit $N$ and $N + 1$.

2.3. *The infinite-volume measure.* The purpose of this paragraph is to extend $\mathcal{P}_\varepsilon$, introduced in Section 2.2, to a law for the whole field $\{\varphi_i\}_{i\in\mathbb{Z}^+}$.

Consider first the regime $\varepsilon \geq \varepsilon_c$, in which case $\tau_k < \infty$ for every $k \in \mathbb{N}$, $\mathcal{P}_\varepsilon$-a.s. We introduce the excursions $\{e_k\}_{k\in\mathbb{N}}$ of the field over the contact set by

$$(2.16) \qquad e_k = \{e_k(i)\}_{0\leq i\leq \tau_k - \tau_{k-1}} := \{\varphi_{\tau_{k-1}+i}\}_{0\leq i\leq \tau_k - \tau_{k-1}}.$$

The variables $e_k$ take values in the space $\bigcup_{m=2}^\infty \mathbb{R}^m$. It is clear that the whole field $\{\varphi_i\}_{i\in\mathbb{Z}^+}$ is in one-to-one correspondence with the process $\{(\tau_k, J_k, e_k)\}_{k\in\mathbb{Z}^+}$. $\mathcal{P}_\varepsilon$ has already been defined as a law for $\{(\tau_k, J_k)\}_{k\in\mathbb{Z}^+}$ [see (2.11)] and we now extend it to a law for $\{\varphi_i\}_{i\in\mathbb{Z}^+}$ in the following way: conditionally on $\{(\tau_k, J_k)\}_{k\in\mathbb{Z}^+}$, we declare that the excursions $\{e_k\}_{k\in\mathbb{N}}$ under $\mathcal{P}_\varepsilon$ are independent, with marginal laws given by

$$\begin{aligned}
e_k \quad &\text{under} \quad \mathcal{P}_\varepsilon(\cdot | \{(\tau_i, J_i)\}_{i\in\mathbb{Z}^+}) \\
&\stackrel{d}{=} (Z_0, \ldots, Z_l) \quad \text{under} \quad \mathbf{P}^{(-a,0)}(\cdot | Z_{l-1} = b, Z_l = 0), \\
&\text{where } l = \tau_k - \tau_{k-1}, a = J_{k-1}, b = J_k.
\end{aligned} \quad (2.17)$$

In words, $e_k$ under $\mathcal{P}_\varepsilon$ is distributed like a bridge of the integrated random walk $\{Z_n\}_n$ of length $l = \tau_k - \tau_{k-1}$, with boundary conditions $Z_{-1} = J_{k-1}$, $Z_0 = 0$, $Z_{l-1} = J_k$ and $Z_l = 0$. Recall in fact that, by (2.7), we have $\mathbf{P}^{(-a,0)} = \mathbf{P}(\cdot | Z_{-1} = a, Z_0 = 0)$, and this is the reason for the minus sign.

Next we consider the regime $\varepsilon < \varepsilon_c$, in which the process $\{\tau_k\}_k$ is $\mathcal{P}_\varepsilon$-a.s. terminating, that is, there is some random index $k^* \in \mathbb{N}$ such that $\tau_k < \infty$ for $k \leq k^*$, while $\tau_{k^*+1} = \infty$. Conditionally on $\{(\tau_k, J_k)\}_{k \in \mathbb{Z}^+}$, the law of the variables $\{e_k\}_{1 \leq k \leq k^*}$ under $\mathcal{P}_\varepsilon$ is still given by (2.17), and to reconstruct the full field $\{\varphi_i\}_{i \in \mathbb{Z}^+}$, it remains to define the law of the last excursion $e_{k^*+1} := \{\varphi_{\tau_{k^*}+i}\}_{0 \leq i < \infty}$, which we do in the following way:

$$e_{k^*+1} \quad \text{under} \quad \mathcal{P}_\varepsilon(\cdot | \{(\tau_i, J_i)\}_{i \in \mathbb{Z}^+}) \stackrel{d}{=} \{Z_i\}_{0 \leq i < \infty} \quad \text{under} \quad \mathbf{P}^{(-J_{k^*}, 0)}(\cdot).$$

This completes the definition of $\mathcal{P}_\varepsilon$ as a law for $\{\varphi_i\}_{i \in \mathbb{Z}^+}$.

Now notice that, conditionally on $\{\ell_N, (\tau_k, J_k)_{k \leq \ell_N}\}$, the excursions $\{e_k\}_{k \leq \ell_N}$ under the pinning model $\mathbb{P}_{\varepsilon,N}$ are independent and their marginal laws are given exactly by (2.17). To see this, it suffices to condition equation (2.5) on $\{J_k\}_{k \leq \ell_N}$, obtaining

$$\mathbb{P}_{\varepsilon,N}(\cdot | \ell_N, (\tau_k, J_k)_{k \leq \ell_N})$$
$$= \mathbf{P}^{(0,0)}((Z_1, \ldots, Z_{N-1}) \in \cdot | Z_{\tau_i} = 0, Z_{\tau_i - 1} = J_k, \ \forall i \leq \ell_N).$$

Then, using the fact that the process $\{Z_n\}_{n \in \mathbb{Z}^+}$ has memory $m = 2$ [see (2.6)], this equation yields easily that the excursions $\{e_k\}_{k \leq \ell_N}$ are indeed conditionally independent and distributed according to (2.17).

These observations have the following important consequence: the basic relation (2.15) can be now extended to hold for the whole field, that is,

$$(2.18) \qquad \mathbb{P}_{\varepsilon,N}(\mathrm{d}\varphi_1, \ldots, \mathrm{d}\varphi_{N-1}) = \mathcal{P}_\varepsilon(\mathrm{d}\varphi_1, \ldots, \mathrm{d}\varphi_{N-1} | \mathcal{A}_N).$$

(Of course, the extension of $\mathcal{P}_\varepsilon$ has been given exactly with this purpose.) Thus, the polymer measure $\mathbb{P}_{\varepsilon,N}$ is nothing but the conditioning of an explicit law $\mathcal{P}_\varepsilon$ with respect to the event $\mathcal{A}_N$. We stress that $\mathcal{P}_\varepsilon$ does not have any dependence on $N$: in this sense, the law $\mathbb{P}_{\varepsilon,N}$ depends on $N$ only through the conditioning on the event $\mathcal{A}_N$. This fact plays a fundamental role in the rest of the paper.

REMARK 2.1. Although the law $\mathcal{P}_\varepsilon$ has been introduced in a somewhat artificial way, it actually has a natural interpretation: it is the *infinite volume limit* of the pinning model, that is, as $N \to \infty$, the law $\mathbb{P}_{\varepsilon,N}$ converges weakly on $\mathbb{R}^{\mathbb{Z}^+}$ to $\mathcal{P}_\varepsilon$. This fact provides another path characterization of the phase transition, because the process $\{\varphi_n\}_{n \in \mathbb{N}}$ under $\mathcal{P}_\varepsilon$ is positive recurrent, null recurrent or transient respectively when $\varepsilon > \varepsilon_c$, $\varepsilon = \varepsilon_c$ or $\varepsilon < \varepsilon_c$. We also note that the field $\{\varphi_i\}_{i \geq 0}$ under the law $\mathcal{P}_\varepsilon$ is not a Markov process, but it rather is a *time-homogeneous* process with finite memory $m = 2$, like $\{Z_n\}_{n \geq 0}$ under $\mathbf{P}$, compare (2.6). Although we do not prove these results, it may be helpful to keep them in mind.

**3. Proof of Theorem 1.4: first part.** In this section we prove a first half of Theorem 1.4, more precisely, (1.14) and the upper bound on $\max_{0 \le k \le N} |\varphi_k|$ in (1.15). Note that these results yield as an immediate corollary the proof of Theorem 1.2 for $\varepsilon \ge \varepsilon_c$ (the case $\varepsilon < \varepsilon_c$ is deferred to Section 4).

The basic tools we use are the description of the pinning law $\mathbb{P}_{\varepsilon,N}$ given in Section 2, that we further develop in Section 3.1 to extract a genuine renewal structure, and a bound based on the Brascamp–Lieb inequality, that we recall in Section 3.2.

3.1. *From Markov renewals to true renewals.* It is useful to observe that, in the framework of Markov renewal theory described in Section 2.2, one can isolate a genuine renewal process. To this purpose, we introduce the (random) set $\chi$ of the adjacent contact points, defined by

(3.1) $$\chi := \{i \in \mathbb{Z}^+ : \varphi_{i-1} = \varphi_i = 0\},$$

and we set by definition $\varphi_{-1} = \varphi_0 = 0$, so that $\chi \ni 0$. We identify $\chi$ with the sequence of random variables $\{\chi_k\}_{k \in \mathbb{Z}^+}$ defined by

(3.2) $$\chi_0 := 0, \qquad \chi_{k+1} := \inf\{i > \chi_k : \varphi_{i-1} = \varphi_i = 0\}, \qquad k \in \mathbb{Z}^+,$$

and we denote by $\iota_N$ the number of adjacent contact points occurring before $N$:

(3.3) $$\iota_N := \#\{\chi \cap [1, N]\} = \sup\{k \in \mathbb{Z}^+ : \chi_k \le N\}.$$

The first observation is that, for every $\varepsilon > 0$, the process $\{\chi_k\}_{k \in \mathbb{Z}^+}$ under the law $\mathcal{P}_\varepsilon$ is a *genuine renewal process*, that is, the increments $\{\chi_{k+1} - \chi_k\}_{k \in \mathbb{Z}^+}$ are independent and identically distributed random variables, taking values in $\mathbb{N} \cup \{\infty\}$, as it is proven in Proposition 5.1 in [6]. Denoting by $q_\varepsilon(n)$ the law of $\chi_1$,

(3.4) $$q_\varepsilon(n) := \mathcal{P}_\varepsilon(\chi_1 = n),$$

it turns out that the properties of $q_\varepsilon(n)$ resemble closely those of $K^\varepsilon_{x,\mathrm{d}y}(n)$, given in Section 2.2. In fact, $q_\varepsilon(\cdot)$ is defective for $\varepsilon < \varepsilon_c$ $[\sum_{n \in \mathbb{N}} q_\varepsilon(n) < 1]$, while it is proper for $\varepsilon \ge \varepsilon_c$ $[\sum_{n \in \mathbb{N}} q_\varepsilon(n) = 1]$. About the asymptotic behavior of $q_\varepsilon(\cdot)$, there exists $\alpha > 0$ such that for every $\varepsilon \in (0, \varepsilon_c + \alpha]$ as $n \to \infty$

(3.5) $$q_\varepsilon(n) = \frac{C_\varepsilon}{n^2} \exp(-\mathrm{F}(\varepsilon) \cdot n)(1 + o(1)),$$

where $C_\varepsilon > 0$; cf. Proposition 7.1 in [6] [which is stated for $\varepsilon \in [\varepsilon_c, \varepsilon_c + \alpha]$, but its proof goes true without changes also for $\varepsilon \in (0, \varepsilon_c)$]. We stress that $\mathrm{F}(\varepsilon) = 0$ for $\varepsilon \le \varepsilon_c$, while $\mathrm{F}(\varepsilon) > 0$ for $\varepsilon > \varepsilon_c$. When $\varepsilon > \varepsilon_c + \alpha$, we content ourselves with the rougher bound

(3.6) $$q_\varepsilon(n) \le C \exp(-\mathrm{G}(\varepsilon) \cdot n) \qquad \forall n \in \mathbb{N},$$

for a suitable $G(\varepsilon) > 0$, which can also be extracted from the proof of Proposition 7.1 in [6] [we have $G(\varepsilon) > F(\varepsilon)$ for large $\varepsilon$].

To summarize, the renewal process $\{\chi_k\}_k$ under $\mathcal{P}_\varepsilon$ is defective with heavy tails in the delocalized regime ($\varepsilon < \varepsilon_c$) and it is proper with heavy tails in the critical regime ($\varepsilon = \varepsilon_c$), while it is proper with exponential tails in the localized regime ($\varepsilon > \varepsilon_c$).

Coming back to the pinning model $\mathbb{P}_{\varepsilon,N}$, by projecting the basic relation (2.15) on the set $\chi$, we obtain that the vector $\{\iota_N, (\chi_k)_{k \le \iota_N}\}$ has the same distribution under $\mathbb{P}_{\varepsilon,N}$ and under $\mathcal{P}_\varepsilon(\cdot|\mathcal{A}_N)$, where we can express the event $\mathcal{A}_N$ in terms of $\chi$, since $\mathcal{A}_N = \{N+1 \in \chi\}$. In words, the adjacent contact points $\{\chi_n\}_n$ under the polymer measure $\mathbb{P}_{\varepsilon,N}$ are distributed like a genuine renewal process conditioned to hit $N+1$.

3.2. *The Brascamp–Lieb inequality.* Let $H:\mathbb{R}^n \to \mathbb{R} \cup \{+\infty\}$ be a function that can be written as

(3.7) $$H(x) = \tfrac{1}{2}A(x) + R(x),$$

where $A(x)$ is a positive definite quadratic form and $R(x)$ is a convex function. Consider the probability laws $\mu_H$ and $\mu_A$ on $\mathbb{R}^n$ defined by

$$\mu_H(\mathrm{d}x) := \frac{e^{-H(x)}}{c_H}\,\mathrm{d}x, \qquad \mu_A(\mathrm{d}x) := \frac{(\det A)^{1/2}}{(2\pi)^{n/2}} e^{-1/2 A(x)}\,\mathrm{d}x,$$

where $\mathrm{d}x$ denotes the Lebesgue measure on $\mathbb{R}^n$ and $c_H$ is the normalizing constant. Of course, $\mu_A$ is a Gaussian law with zero mean and with $A^{-1}$ as covariance matrix.

We denote by $E_H$ and $E_A$ respectively the expectation with respect to $\mu_H$ and $\mu_A$. The Brascamp–Lieb inequality reads as follows (cf. [5], Corollary 6): for any convex function $\Gamma:\mathbb{R} \to \mathbb{R}$ and for all $a \in \mathbb{R}^n$, such that $E_A[\Gamma(a \cdot x)] < \infty$, we have

(3.8) $$E_H[\Gamma(a \cdot x - E_H(a \cdot x))] \le E_A[\Gamma(a \cdot x)],$$

where $a \cdot x$ denotes the standard scalar product on $\mathbb{R}^n$.

A useful observation is that (3.8) still holds true if we condition $\mu_H$ through linear constraints. More precisely, given $m \le n$ and $b_i \in \mathbb{R}^n$, $c_i \in \mathbb{R}$ for $1 \le i \le m$, we set

$$\mu_H^*(\mathrm{d}x) := \mu_H(\mathrm{d}x|b_i \cdot x = c_i,\ \forall i \le m).$$

We assume that the set of solutions of the linear system $\{b_i \cdot x = c_i,\ \forall i \le m\}$ has nonempty intersection with the support of $\mu_H$ and that $x \mapsto e^{-H(x)}$ is continuous on the whole $\mathbb{R}^n$, so that there is no problem in defining the conditional measure $\mu_H^*$. Let us proceed through an approximation argument: for $k \in \mathbb{N}$ we set

$$H_k^*(x) := H(x) + k \sum_{i=1}^m (b_i \cdot x - c_i)^2, \qquad \mu_{H_k}^*(\mathrm{d}x) := \frac{e^{-H_k(x)}}{c_{H_k}^*}\,\mathrm{d}x,$$

where $c^*_{H_k}$ is the normalizing constant that makes $\mu^*_{H_k}$ a probability. Since we have added convex terms, $H^*_k(x)$ is still of the form (3.7), with the same $A(x)$, hence, equation (3.8) holds true with $E_H$ replaced by $E^*_{H_k}$. However, it is easy to realize that $\mu^*_{H_k}$ converges weakly to $\mu^*_H$ as $k \to \infty$, hence, (3.8) holds true also for $E^*_H$, that is,

$$(3.9) \qquad E^*_H[\Gamma(a \cdot x - E^*_H(a \cdot x))] \leq E_A[\Gamma(a \cdot x)].$$

3.3. *A preliminary bound.* Before passing to the proof of Theorem 1.4, we derive a useful bound based on the Brascamp–Lieb inequality. We recall that, by the uniform strict convexity assumption on the potential, we can write $V(t) = \frac{\gamma}{2}t^2 + r(t)$, where $\gamma > 0$ [cf. (1.1)] and $r(\cdot) : \mathbb{R} \to \mathbb{R} \cup \{+\infty\}$ is a convex function; see Section 1.1.

By (2.4), the law of the vector $(Z_1, \ldots, Z_n)$ under $\mathbf{P}^{(0,0)}$ has the form $\mu_H(\mathrm{d}x) = e^{-H(x)}\,\mathrm{d}x$, $x \in \mathbb{R}^n$, where $H(x) = \frac{1}{2}A(x) + R(x)$ with

$$(3.10) \quad A(x_1, \ldots, x_n) = \gamma \cdot \left( (x_1)^2 + (x_2 - 2x_1)^2 + \sum_{i=2}^{n-1}(x_{i+1} + x_{i-1} - 2x_i)^2 \right),$$

$$R(x_1, \ldots, x_n) = r(x_1) + r(x_2 - 2x_1) + \sum_{i=2}^{n-1} r(x_{i+1} + x_{i-1} - 2x_i).$$

Since $r(\cdot)$ is convex on $\mathbb{R}$, $R(\cdot)$ is convex on $\mathbb{R}^n$ and, therefore, we are in the Brascamp–Lieb framework described in Section 3.2. Fix arbitrarily $m \leq n$ and $t_1, \ldots, t_m \in \{1, \ldots, n\}$ and consider $\mu^*_H(\mathrm{d}x) = \mu_H(\mathrm{d}x|x_{t_1} = 0, \ldots, x_{t_m} = 0)$. Applying (3.9) with $\Gamma(x) = e^{\lambda x_k}$, for $\lambda \in \mathbb{R}$ and $k \in \{1, \ldots, n\}$, and noting that $E^*_H(x_k) = 0$ by symmetry, we obtain

$$\mathbf{E}^{(0,0)}(e^{\lambda Z_k}|Z_{t_1} = 0, \ldots, Z_{t_m} = 0)$$
$$= E^*_H(e^{\lambda x_k}) \leq E_A(e^{\lambda x_k}) = \exp\left(\frac{\lambda^2}{2\gamma} \cdot \frac{k(k+1)(2k+1)}{6}\right),$$

where the last equality is the result of a straightforward Gaussian computation, because in this context $\mu_A$ is just the law of the integral of a random walk with Gaussian steps $\sim \mathcal{N}(0, \gamma^{-1})$ [cf. (2.2) and (2.4)]. Applying Markov's inequality and optimizing over $\lambda$ yields for $s \in \mathbb{R}^+$

$$(3.11) \qquad \mathbf{P}^{(0,0)}(|Z_k| > s|Z_{t_1} = 0, \ldots, Z_{t_m} = 0) \leq 2\exp\left(-\frac{\gamma}{k^3}\frac{s^2}{6}\right).$$

The crucial aspect is that this bound is uniform over the choices of the points $t_i$ (that do not appear in the r.h.s.). In a sense, this is no surprise, because conditioning on $Z_{t_i} = 0$ should decrease the probability of the event $\{|Z_k| > t\}$.

3.4. *Proof of Theorem 1.4: upper bounds.* Recalling the basic relation (2.5), for $m \leq N-1$ and $t_1, \ldots, t_m \in \{1, \ldots, N-1\}$ we have

$$
\begin{aligned}
&\mathbb{P}_{\varepsilon,N}(|\varphi_k| > s \mid \ell_N = m, \tau_i = t_i, \ \forall 1 \leq i \leq m) \\
&= \mathbf{P}^{(0,0)}(|Z_k| > s \mid Z_j = 0, \ \forall j \in \{t_1, t_2, \ldots, t_m\} \cup \{N, N+1\}).
\end{aligned}
\tag{3.12}
$$

We now observe that the process $\{Z_n\}_n$ under $\mathbf{P}^{(0,0)}$ is a process with finite memory $m = 2$ [see (2.6)], hence, its excursions between adjacent zeros are independent. For this reason, we identify the adjacent zeros that are close to $k$, in the following way: we first set, for convenience,

$$t_{-1} := -1, \qquad t_0 := 0, \qquad t_{m+1} := N, \qquad t_{m+2} := N+1$$

and we define

$$\mathbf{l} := \max\{i \geq 0 : t_i \leq k \text{ and } t_i - t_{i-1} = 1\},$$
$$\mathbf{r} := \min\{i \geq 0 : t_i > k \text{ and } t_i - t_{i-1} = 1\}.$$

In words, $t_{\mathbf{l}}$ (resp. $t_{\mathbf{r}}$) is the closest adjacent zero at the left (resp. at the right) of $k$. Note that $0 \leq \mathbf{l} < \mathbf{r} \leq m+1$. Then the above mentioned finite memory property yields

$$
\begin{aligned}
&\mathbf{P}^{(0,0)}(|Z_k| > s \mid Z_j = 0, \ \forall j \in \{t_1, \ldots, t_m\} \cup \{N, N+1\}) \\
&= \mathbf{P}^{(0,0)}(|Z_k| > s \mid Z_j = 0, \ \forall j \in \{t_{\mathbf{l}-1}, t_{\mathbf{l}}, \ldots, t_{\mathbf{r}}\}) \\
&= \mathbf{P}^{(0,0)}(|Z_{k-t_{\mathbf{l}}}| > s \mid Z_j = 0, \ \forall j \in \{t_{\mathbf{l}+1} - t_{\mathbf{l}}, t_{\mathbf{l}+2} - t_{\mathbf{l}}, \ldots, t_{\mathbf{r}} - t_{\mathbf{l}}\}),
\end{aligned}
\tag{3.13}
$$

where the second inequality follows by time homogeneity. Putting together (3.12), (3.13) and (3.11), we get

$$
\mathbb{P}_{\varepsilon,N}(|\varphi_k| > s \mid \ell_N = m, \tau_i = t_i, \ \forall 1 \leq i \leq m) \leq 2\exp\left(-\frac{\gamma}{(k-t_{\mathbf{l}})^3}\frac{s^2}{6}\right).
$$
(3.14)

We denote by $\delta_N$ the maximal gap in the adjacent contact set $\chi$ until $N$, that is,

$$\delta_N := \max\{\chi_k - \chi_{k-1} : 0 < k \leq \iota_N\}, \tag{3.15}$$

where the variable $\iota_N$ was introduced in (3.3). Then the bound (3.14) yields finally

$$\mathbb{P}_{\varepsilon,N}(|\varphi_k| > s \mid \tau \cap (0,N)) \leq 2\exp\left(-\frac{\gamma}{6(\delta_N)^3}s^2\right). \tag{3.16}$$

This is the key estimate to prove the upper bounds in (1.14) and (1.15). In fact, the inclusion bound yields

$$\mathbb{P}_{\varepsilon,N}\left(\max_{k=1,\ldots,N}|\varphi_k| > s \mid \tau \cap (0,N)\right) \leq 2N\exp\left(-\frac{\gamma}{6(\delta_N)^3}s^2\right). \tag{3.17}$$

It is now clear the importance of studying the asymptotic behavior of the variable $\delta_N$.

We start considering the critical regime ($\varepsilon = \varepsilon_c$). As we prove in Appendix A.1, there exists a positive constant $c_1$ and a sequence $(a_n)_n$ such that, for all $N \geq 3$ and $t \in [1, \infty)$,

$$(3.18) \quad \mathbb{P}_{\varepsilon_c, N}\left(\delta_N \geq t \frac{N}{\log N}\right) \leq \frac{c_1}{t} + a_N \quad \text{with } a_N \to 0 \text{ as } N \to \infty.$$

Combining this relation with (3.17), we get

$$\mathbb{P}_{\varepsilon_c, N}\left(\max_{k=1,\ldots,N} |\varphi_k| > s\right)$$

$$\leq \mathbb{P}_{\varepsilon_c, N}\left(\max_{k=1,\ldots,N} |\varphi_k| > s, \delta_N < t \frac{N}{\log N}\right) + \frac{c_1}{t} + a_N$$

$$\leq 2N \mathbb{E}_{\varepsilon_c, N}\left[\exp\left(-\frac{\gamma}{6(\delta_N)^3} s^2\right) \mathbf{1}_{\{\delta_N < tN/\log N\}}\right] + \frac{c_1}{t} + a_N$$

$$\leq 2N \exp\left(-\frac{\gamma}{6t^3} \frac{(\log N)^3}{N^3} s^2\right) + \frac{c_1}{t} + a_N$$

and, setting $s = KN^{3/2}/\log N$ and $t = (\frac{\gamma}{12})^{1/3} K^{2/3}$, we finally obtain

$$\mathbb{P}_{\varepsilon_c, N}\left(\max_{k=1,\ldots,N} |\varphi_k| > K \frac{N^{3/2}}{\log N}\right) \leq \frac{2}{N} + \frac{c_1(12/\gamma)^{1/3}}{K^{2/3}} + a_N.$$

Since $a_N \to 0$ as $N \to \infty$ [see (3.18)], the upper bound in equation (1.15) is proven.

Then we consider the localized regime ($\varepsilon > \varepsilon_c$). As we prove in Appendix A.3, there exists a positive constant $c_2$ such that

$$(3.19) \qquad \mathbb{P}_{\varepsilon, N}(\delta_N \geq c_2 \log N) \longrightarrow 0 \quad \text{as } N \to \infty.$$

Then, in analogy with the preceding lines, we combine this relation with (3.17), getting

$$\mathbb{P}_{\varepsilon, N}\left(\max_{k=1,\ldots,N} |\varphi_k| > s\right) \leq \mathbb{P}_{\varepsilon, N}\left(\max_{k=1,\ldots,N} |\varphi_k| > s, \delta_N < c_2 \log N\right) + o(1)$$

$$\leq 2N \mathbb{E}_{\varepsilon, N}\left[\exp\left(-\frac{\gamma}{6(\delta_N)^3} s^2\right) \mathbf{1}_{\{\delta_N < c_2 \log N\}}\right] + o(1)$$

$$\leq 2N \exp\left(-\frac{\gamma}{6(c_2)^3} \frac{s^2}{(\log N)^3}\right) + o(1).$$

Setting $s = K(\log N)^2$, for $K$ sufficiently large we obtain

$$\mathbb{P}_{\varepsilon_c, N}\left(\max_{k=1,\ldots,N} |\varphi_k| > K(\log N)^2\right) \leq 2N^{-\gamma K^2/(6(c_2)^3)+1} + o(1) \longrightarrow 0$$

as $N \to \infty$,

hence, also (1.14) is proven.

**4. Proof of Theorems 1.3 and 1.2.** In this section we focus on the delocalized regime $\varepsilon < \varepsilon_c$, proving Theorem 1.3 and the corresponding part of Theorem 1.2. We recall that the proof of Theorem 1.2 for $\varepsilon \geq \varepsilon_c$ follows immediately from the upper bound on $\max_{0 \leq k \leq N} |\varphi_k|$ given by relations (1.14) and (1.15), that have already been proven in Section 3.

4.1. *The free case $\varepsilon = 0$.* We start proving Theorem 1.2 in the case $\varepsilon = 0$, when there is no interaction between the field and the defect line. The main ingredient is the random walk interpretation outlined in Section 2.1. We recall from Section 1.3 that $\{B_t\}_{t \in [0,1]}$ denotes a standard Brownian motion on $\mathbb{R}$ and $I_t = \int_0^t B_s \, ds$ denotes its integral, while $(\widehat{B}_t, \widehat{I}_t)$ denotes $(B_t, I_t)$ conditionally on $(B_1, I_1) = (0, 0)$; see (1.11).

We first state a local limit theorem for the process $\{Z_n\}_{n \in \mathbb{N}}$, proven in Proposition 2.3 of [6]. We note that the vector $(Y_n, Z_n) = (Z_n - Z_{n-1}, Z_n)$ has an absolutely continuous law under $\mathbf{P}^{(a,b)}$ for $n \geq 2$, and we introduce its density

$$\varphi_n^{(a,b)}(y, z) := \frac{\mathbf{P}^{(a,b)}((Y_n, Z_n) \in (\mathrm{d}y, \mathrm{d}z))}{\mathrm{d}y \, \mathrm{d}z}.$$

Notice that $\varphi_n^{(a,b)}(y, z) = \varphi_n^{(0,0)}(y - a, z - b - na)$ by (2.3), hence, we can focus on $\varphi_n^{(0,0)}$. The local limit theorem reads as follows:

$$(4.1) \quad \sup_{(y,z) \in \mathbb{R}^2} |\sigma^2 n^2 \varphi_n^{(0,0)}(\sigma \sqrt{n} y, \sigma n^{3/2} z) - g(y, z)| \longrightarrow 0 \qquad (n \to \infty),$$

where $g(y, z) := \frac{\sqrt{3}}{\pi} \exp(-2y^2 - 6z^2 + 6yz)$ is the law of the Gaussian vector $(B_1, I_1)$.

We are ready to prove a somewhat general invariance principle, from which Theorem 1.2 for $\varepsilon = 0$ follows as a corollary, because $\mathbb{P}_{0,N}$ coincides with the law of the integrated random walk $(Z_1, \ldots, Z_{N-1})$ under $\mathbf{P}^{(0,0)}(\cdot | Y_{N+1} = 0, Z_{N+1} = 0)$; compare Section 2.1. For notational convenience, we simply denote by $Z_{\langle Nt \rangle}$ and $Y_{\langle Nt \rangle}$ the linear interpolation of the processes.

PROPOSITION 4.1. *Uniformly for $a, c$ in compact sets of $\mathbb{R}$, we have, as $N \to \infty$,*

$$(4.2) \quad \left\{ \left( \frac{Y_{\langle Nt \rangle}}{\sigma \sqrt{N}}, \frac{Z_{\langle Nt \rangle}}{\sigma N^{3/2}} \right) \right\}_{t \in [0,1]}$$

$$under \quad \mathbf{P}^{(-a,0)}(\cdot | (Y_N, Z_N) = (-c, 0)) \xrightarrow{d} \{(\widehat{B}_t, \widehat{I}_t)\}_{t \in [0,1]},$$

*where $\xrightarrow{d}$ denotes convergence in distribution on $C([0,1]) \times C([0,1])$.*

PROOF. We start noting that the process $\{(\frac{Y_{\langle Nt\rangle}}{\sigma\sqrt{N}}, \frac{Z_{\langle Nt\rangle}}{\sigma N^{3/2}})\}_{t\in[0,1]}$ under the unconditioned law $\mathbf{P}^{(-a,0)}(\cdot)$ converges in distribution as $N \to \infty$ toward $\{(B_t, I_t)\}_{t\in[0,1]}$, uniformly for $a$ in compact sets of $\mathbb{R}$. This is an easy consequence of Donsker's Invariance Principle and the Continuous Mapping Theorem, because $\{Y_n\}_{n\in\mathbb{Z}^+}$ under $\mathbf{P}^{(-a,0)}$ is a zero-mean, finite-variance real random walk starting at $-a$ and, moreover,

$$\frac{Z_{\langle Nt\rangle}}{\sigma N^{3/2}} = \int_0^t \frac{(Y_{\lceil Ns\rceil} + a)}{\sigma\sqrt{N}}\,\mathrm{d}s$$

(we recall that $(\xi_t)_t \mapsto \int_0^t \xi_s\,\mathrm{d}s$ is a continuous functional on $D([0,1])$).

Next it is convenient to restrict the parameter $t$ to $[0, 1-\eta]$, where $\eta > 0$ is fixed. Since $\{(Y_n, Z_n)\}_{n\in\mathbb{Z}^+}$ is a Markov process, the law of the process $\{(\frac{Y_{\langle Nt\rangle}}{\sigma\sqrt{N}}, \frac{Z_{\langle Nt\rangle}}{\sigma N^{3/2}})\}_{t\in[0,1-\eta]}$ under $\mathbf{P}^{(-a,0)}(\cdot|(Y_N, Z_N) = (-c, 0))$ is absolutely continuous w.r.t. the law of the same process under $\mathbf{P}^{(-a,0)}(\cdot)$, with Radon–Nikodym derivative $f_N^{(\eta)}$ given by

$$f_N^{(\eta)}((y_t, z_t)_{t\in[0,1-\eta]}) = f_N^{(\eta)}(y_{1-\eta}, z_{1-\eta}) = \frac{\varphi_{\lfloor \eta N\rfloor}^{(\sigma\sqrt{N}y_{1-\eta}, \sigma N^{3/2}z_{1-\eta})}(-c, 0)}{\varphi_N^{(-a,0)}(-c, 0)}.$$

The local limit theorem (4.1) yields the uniform convergence on compact sets of the function $f_N^{(\eta)}$ as $N \to \infty$ toward an explicit limit function $f^{(\eta)}$, uniformly for $a, c$ in compact sets, and one checks directly that $f^{(\eta)}$ is indeed the Radon–Nikodym derivative of the law of $\{(\widehat{B}_t, \widehat{I}_t)\}_{t\in[0,1-\eta]}$ w.r.t. the law of $\{(B_t, I_t)\}_{t\in[0,1-\eta]}$. This shows that equation (4.2) holds when $t$ is restricted to $[0, 1-\eta]$. Since this is true for every $\eta > 0$, the proof is completed with a standard tightness argument. $\square$

4.2. *Proof of Theorem 1.3.* In this paragraph we focus on the regime $0 < \varepsilon < \varepsilon_c$. We start proving a slightly stronger version of equation (1.13). We denote by $l_N$ (resp. $r_N$) the index of the last point in the contact set before $N/2$ (resp. after $N/2$), that is,

(4.3) $\quad l_N := \max\{i \geq 0 : \tau_i \leq N/2\}, \qquad r_N := \min\{i \geq 0 : \tau_i > N/2\} = l_N + 1.$

Equation (1.13) says that both $\tau_{l_N}$ and $N - \tau_{r_N}$ are $O(1)$. It turns out that also $|J_{l_N}|$ and $|J_{r_N}|$ are $O(1)$. More precisely, we are going to prove that

(4.4) $\quad \lim_{L\to+\infty} \liminf_{N\to+\infty} \mathbb{P}_{\varepsilon,N}(\tau_{l_N} \leq L, \tau_{r_N} \geq N - L, |J_{l_N}| \leq L, |J_{r_N}| \leq L) = 1.$

As a matter of fact, it is not difficult to further strengthen this relation, by showing that also $\max_{0\leq i\leq \tau_{l_N}} |\varphi_i|$ and $\max_{\tau_{r_N}\leq i\leq N} |\varphi_i|$ are $O(1)$, but we omit the details for conciseness.

The proof of (4.4) is based on relation (2.18) [or, more directly, on (2.15)]. Recalling the definition (2.11) of the transition kernel $\mathsf{K}^\varepsilon_{x,\mathrm{d}y}(n)$, we introduce the associated renewal kernel $\mathsf{U}^\varepsilon_{x,\mathrm{d}y}(n)$ by

$$(4.5) \qquad \mathsf{U}^\varepsilon_{x,\mathrm{d}y}(n) := \sum_{k=0}^\infty (\mathsf{K}^\varepsilon)^{*k}_{x,\mathrm{d}y}(n) = \mathcal{P}_\varepsilon(n \in \tau, \varphi_{n-1} \in \mathrm{d}y | J_0 = x),$$

where $(\mathsf{K}^\varepsilon)^{*k}$ denotes the $k$-fold convolution of the kernel $\mathsf{K}^\varepsilon$ with itself: by definition, $(\mathsf{K}^\varepsilon)^{*0}_{x,\mathrm{d}y}(n) := \delta_x(\mathrm{d}y)\mathbf{1}_{\{n=0\}}$ and, given the kernels $\mathsf{A}_{x,\mathrm{d}y}(n), \mathsf{B}_{x,\mathrm{d}y}(n)$, we set

$$(\mathsf{A} * \mathsf{B})_{x,\mathrm{d}y}(n) := \sum_{m=0}^n \int_{z\in\mathbb{R}} \mathsf{A}_{x,\mathrm{d}z}(m) \cdot \mathsf{B}_{z,\mathrm{d}y}(n-m).$$

In particular, $\mathcal{P}_\varepsilon(\mathcal{A}_N) = U_{0,\{0\}}(N+1)$. With this notation, by (2.18), we can write

$$\mathbb{P}_{\varepsilon,N}(\tau_{l_N} \leq L, \tau_{r_N} \geq N-L, |J_{l_N}| \leq L, |J_{r_N}| \leq L)$$

$$(4.6) \qquad = \frac{1}{U^\varepsilon_{0,\{0\}}(N+1)}$$

$$\times \sum_{a,b=0}^L \int_{x,y\in[-L,L]} \mathsf{U}^\varepsilon_{0,\mathrm{d}x}(a) \cdot \mathsf{K}^\varepsilon_{x,\mathrm{d}y}(N+1-a-b) \cdot \mathsf{U}^\varepsilon_{y,\{0\}}(b).$$

By (2.12) and (4.1), it follows that, for bounded $x,y$ and as $n \to \infty$,

$$(4.7) \qquad \mathsf{K}^\varepsilon_{x,\mathrm{d}y}(n) \sim \frac{L^\varepsilon_{x,\mathrm{d}y}}{n^2}, \qquad \text{where } L^\varepsilon_{x,\mathrm{d}y} := \frac{6\varepsilon}{\pi} \frac{v_\varepsilon(y)}{v_\varepsilon(x)} \mathrm{d}y.$$

To determine the asymptotic behavior of $U_{0,\{0\}}(N+1)$ as $N \to \infty$, we apply the Markov Renewal Theorem given by equation (7.9) in [6], Section 7.2 (it is easily checked that all the assumptions are verified). We set

$$(4.8) \qquad B^\varepsilon_{x,\mathrm{d}y} := \sum_{n\in\mathbb{N}} \mathsf{K}^\varepsilon_{x,\mathrm{d}y}(n), \qquad (1-B^\varepsilon)^{-1}_{x,\mathrm{d}y} := \sum_{k=0}^\infty (B^\varepsilon)^{\circ k}_{x,\mathrm{d}y},$$

where $(B^\varepsilon)^{\circ k}$ denotes the $k$-fold composition of the kernel $B^\varepsilon$ with itself: by definition, $(B^\varepsilon)^{\circ 0}_{x,\mathrm{d}y} := \delta_x(\mathrm{d}y)$ and $(A \circ B)_{x,\mathrm{d}y} := \int_{z\in\mathbb{R}} A_{x,\mathrm{d}z} B_{z,\mathrm{d}y}$. Then by equation (7.9) in [6], we can write, as $n \to \infty$,

$$U^\varepsilon_{0,\{0\}}(n) \sim \frac{((1-B^\varepsilon)^{-1} \circ L \circ (1-B^\varepsilon)^{-1})_{0,\{0\}}}{n^2}$$

and, therefore,

$$\lim_{N\to\infty} \mathbb{P}_{\varepsilon,N}(\tau_{l_N} \leq L, \tau_{r_N} \geq N-L, |J_{l_N}| \leq L, |J_{r_N}| \leq L)$$

$$= \frac{\int_{x,y\in[-L,L]} (\sum_{a=0}^L \mathsf{U}^\varepsilon_{0,\mathrm{d}x}(a)) L^\varepsilon_{x,\mathrm{d}y} (\sum_{b=0}^L \mathsf{U}^\varepsilon_{y,\{0\}}(b))}{((1-B^\varepsilon)^{-1} \circ L^\varepsilon \circ (1-B^\varepsilon)^{-1})_{0,\{0\}}}.$$

Since by definition $\sum_{n\in\mathbb{N}} \mathsf{U}^\varepsilon_{x,\mathrm{d}y}(n) = B^\varepsilon_{x,\mathrm{d}y}$, letting $L \to \infty$, the r.h.s. of the last relation converges to 1 and equation (4.4) is proven.

4.3. *Proof of Theorem 1.2.* The proof of Theorem 1.2 for $0 < \varepsilon < \varepsilon_c$ follows by putting together the results proven so far. For conciseness, we just sketch the main arguments and leave the details to the reader.

It is convenient to split the field $\{\varphi_i\}_{0 \le i \le N}$ in three parts: the beginning $\{\varphi_i\}_{0 \le i \le \tau_{l_N}}$, the bulk $\{\varphi_i\}_{\tau_{l_N} \le i \le \tau_{r_N}}$ and the end $\{\varphi_i\}_{\tau_{r_N} \le i \le N}$, where we recall that the indexes $l_N, r_N$ have been introduced in (4.3). By (4.6), both $\tau_{l_N}$ and $N - \tau_{r_N}$ are $O(1)$ as $N \to \infty$. Furthermore, as we already mention, one can show that also $\max_{0 \le i \le \tau_{l_N}} |\varphi_i|$ and $\max_{\tau_{r_N} \le i \le N} |\varphi_i|$ are $O(1)$ as $N \to \infty$. Therefore, both the beginning and the end of the field are irrelevant for the scaling limit [remember the definition (1.10) of the rescaled field $\widehat{\varphi}_N(t)$] and it suffices to focus on the bulk.

We recall that the polymer measure $\mathbb{P}_{\varepsilon,N}$ coincides with the law $\mathcal{P}_\varepsilon$ conditioned on $\mathcal{A}_N$; compare (2.18). In particular, by the construction of $\mathcal{P}_\varepsilon$ given in Sections 2.2 and 2.3, it follows that if we fix $\tau_{l_N} = m$, $\varphi_{m-1} = a$, $\tau_{r_N} = N - n$, $\varphi_{N-n-1} = c$ (of course, $\varphi_m = \varphi_{N-n} = 0$), the bulk $\{\varphi_i\}_{m \le i \le N-n}$ under $\mathbb{P}_{\varepsilon,N}$ is distributed like the process $\{Z_j\}_{0 \le j \le N-n-m}$ under $\mathbf{P}^{(-a,0)}(\cdot|(Y_{N-n-m}, Z_{N-n-m}) = (-c, 0))$. Since all the parameters $m, n, a, c$ are $O(1)$ by (4.4), we can apply Proposition 4.1 and Theorem 1.2 is proven.

**5. Proof of Theorem 1.4: second part.** In this section we complete the proof of Theorem 1.4, by showing that also the lower bound on $\max_{0 \le k \le N} |\varphi_k|$ in (1.15) holds true.

The first basic ingredient, that we prove in Appendix A.2, is a lower bound counterpart of equation (3.18):

$$(5.1) \qquad \lim_{t \to 0^+} \liminf_{N \to \infty} \mathbb{P}_{\varepsilon_c,N}\left(\delta_N \ge t \frac{N}{\log N}\right) = 1.$$

The second ingredient is given by the following lemma, proven in Section 5.1, that will be used also in the proof of Theorem 1.6. Recall the definition (1.10) of the rescaled field $\{\widehat{\varphi}_N(t)\}_{t \in [0,1]}$.

LEMMA 5.1. *Under the conditional law $\mathcal{P}_{\varepsilon_c}(\cdot|\chi_1 = N+1)$, the process $\{\widehat{\varphi}_N(t)\}_{t \in [0,1]}$ converges in distribution on $C([0,1])$ as $N \to \infty$ toward the process $\{\widehat{I}_t\}_{t \in [0,1]}$.*

The idea to complete the proof of Theorem 1.4 is now quite simple. We first notice that, given a gap $(\chi_k, \chi_{k+1})$ in the set $\chi$ of width $m = \chi_{k+1} - \chi_k$, the law of the field inside this gap is nothing but $\mathcal{P}_{\varepsilon_c}(\cdot|\chi_1 = m)$. In particular, by Lemma 5.1, the scaling behavior of the field in this gap is of order $m^{3/2}$.

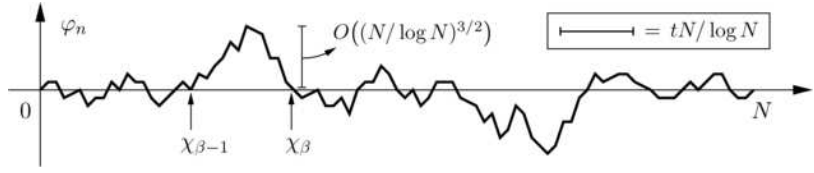

Fig. 2. *A typical trajectory of the field $\{\varphi_n\}_{0\leq n\leq N}$ under the critical law $\mathbb{P}_{\varepsilon_c,N}$. The variables $\chi_{\beta-1}$ and $\chi_\beta$ are the extremities of the first large gap in the set $\chi$ of adjacent contact points; compare (5.4). For simplicity, the distinction between simple and adjacent contact points (i.e., between the sets $\tau$ and $\chi$) is not evidenced in the picture.*

By (5.1), the width of the largest gap in the set $\chi$ before $N$ is of order $\approx N/\log N$, hence, inside this gap the field scales like $(N/\log N)^{3/2}$, from which the lower bound in (1.15) follows. Let us now make these considerations precise (it may be helpful to look at Figure 2).

For $m \in \mathbb{N}$ and $s \in \mathbb{R}^+$ we introduce the event $A_{m,s} := \{\max_{0\leq k\leq m} |\varphi_k| \geq sm^{3/2}\}$, and we note that, by Lemma 5.1, we have

$$(5.2) \qquad \lim_{s\to 0^+} \liminf_{m\to\infty} \mathcal{P}_{\varepsilon_c}(A_{m,s}|\chi_1 = m) = 1.$$

By (5.1), for every $\eta > 0$ we can fix $t > 0$ and $N_0 \in \mathbb{N}$ such that, for all $N \geq N_0$,

$$(5.3) \qquad \mathbb{P}_{\varepsilon_c,N}\left(\delta_N \geq t\frac{N}{\log N}\right) \geq 1 - \eta.$$

We denote by $\beta$ the index of the first long gap in the set $\chi$ (cf. Figure 2):

$$(5.4) \qquad \beta := \inf\left\{i \geq 1 : \chi_i - \chi_{i-1} \geq t\frac{N}{\log N}\right\}.$$

The law of the field inside the gap admits the following explicit description, that follows from relation (2.18): for all $a, b \in \mathbb{N}$ with $0 \leq a \leq b \leq N$ and $b - a \geq tN/\log N$,

$$(5.5) \qquad \begin{aligned} \mathbb{P}_{\varepsilon_c,N}(\{\varphi_i\}_{a\leq i\leq b} \in \cdot |\chi_{\beta-1} = a, \chi_\beta = b) \\ = \mathcal{P}_{\varepsilon_c}(\{\varphi_i\}_{0\leq i\leq b-a} \in \cdot |\chi_1 = b - a). \end{aligned}$$

Observing that $\{\delta_N \geq tN/\log N\} = \{\chi_\beta \leq N\}$ and applying the inclusion bound, we get

$$\mathbb{P}_{\varepsilon_c,N}\left(\max_{0\leq k\leq N} |\varphi_k| \geq \frac{1}{K}\frac{N^{3/2}}{(\log N)^{3/2}}\right)$$
$$\geq \mathbb{P}_{\varepsilon_c,N}\left(\max_{0\leq k\leq N} |\varphi_k| \geq \frac{1}{K}\frac{N^{3/2}}{(\log N)^{3/2}}, \chi_\beta \leq N\right)$$

$$\geq \mathbb{P}_{\varepsilon_c,N}\left(\max_{\chi_{\beta-1}\leq k\leq \chi_\beta}|\varphi_k|\geq \frac{1}{K}\frac{N^{3/2}}{(\log N)^{3/2}}, \chi_\beta \leq N\right)$$

$$= \sum_{\substack{0\leq a\leq b\leq N \\ b-a\geq tN/\log N}} \mathbb{P}_{\varepsilon_c,N}\left(\max_{a\leq k\leq b}|\varphi_k|\geq \frac{1}{K}\frac{N^{3/2}}{(\log N)^{3/2}}, \chi_{\beta-1}=a, \chi_\beta=b\right).$$

Combining this relation with (5.5) and recalling the definition of $A_{m,s}$ yields

$$\mathbb{P}_{\varepsilon_c,N}\left(\max_{0\leq k\leq N}|\varphi_k|\geq \frac{1}{K}\frac{N^{3/2}}{(\log N)^{3/2}}\right)$$
$$\geq \sum_{\substack{0\leq a\leq b\leq N \\ b-a\geq tN/\log N}} \mathcal{P}_{\varepsilon_c}(A_{b-a,1/KN^{3/2}/(\log N)^{3/2}\cdot 1/(b-a)^{3/2}}|\chi_1=b-a)$$
$$\times \mathbb{P}_{\varepsilon_c,N}(\chi_{\beta-1}=a, \chi_\beta=b).$$

Now observe that in the range of summation $\frac{1}{K}\frac{N^{3/2}}{(\log N)^{3/2}} \cdot \frac{1}{(b-a)^{3/2}} \leq \frac{1}{Kt^{3/2}}$ and that the event $A_{m,s}$ is decreasing in $s$. Since $t>0$ is fixed, it follows from (5.2) that for $K$ and $N$ sufficiently large, when $b-a\geq tN/\log N$, we have

$$\mathcal{P}_{\varepsilon_c}(A_{b-a,1/KN^{3/2}/(\log N)^{3/2}\cdot 1/(b-a)^{3/2}}|\chi_1=b-a)$$
$$\geq \mathcal{P}_{\varepsilon_c}(A_{b-a,1/(Kt^{3/2})}|\chi_1=b-a) \geq 1-\eta.$$

Therefore, for the same $K$ and $N$ we get

$$\mathbb{P}_{\varepsilon_c,N}\left(\max_{0\leq k\leq N}|\varphi_k|\geq \frac{1}{K}\frac{N^{3/2}}{(\log N)^{3/2}}\right)$$
$$\geq (1-\eta)\sum_{\substack{0\leq a\leq b\leq N \\ b-a\geq tN/\log N}}\mathbb{P}_{\varepsilon_c,N}(\chi_{\beta-1}=a, \chi_\beta=b)$$
$$= (1-\eta)\mathbb{P}_{\varepsilon_c,N}(\chi_\beta \leq N) \geq (1-\eta)^2,$$

where the last inequality is just (5.3). Since $\eta>0$ was arbitrary, the proof of the lower bound in (1.15) is completed.

5.1. *Proof of Lemma 5.1.* Arguing as in Section 4.3, it suffices to show that under the law $\mathcal{P}_{\varepsilon_c}(\cdot|\chi_1=N+1)$ the contact set is concentrated near the boundary points, and the invariance principle will follow from Proposition 4.1. Recalling the definition (4.3) of the indexes $l_N$ and $r_N$, we prove that

$$\lim_{L\to+\infty}\liminf_{N\to+\infty}\mathcal{P}_{\varepsilon_c}(\tau_{l_N}\leq L, \tau_{r_N}\geq N-L,$$
(5.6)
$$|J_{l_N}|\leq L, |J_{r_N}|\leq L|\chi_1=N+1)=1.$$

Some notation first. We set $\widehat{\mathsf{K}}^\varepsilon_{x,\mathrm{d}y}(n) := \mathsf{K}^\varepsilon_{x,\mathrm{d}y}(n)\mathbf{1}_{\{y\ne 0\}} = \mathsf{K}^\varepsilon_{x,\mathrm{d}y}(n)\mathbf{1}_{\{n\ge 2\}}$ [cf. (2.12)], that gives the law of the jumps occurring before $\chi_1$, and we introduce the corresponding renewal kernel

$$\widehat{\mathsf{U}}^\varepsilon_{x,\mathrm{d}y}(n) := \sum_{k=0}^\infty (\widehat{\mathsf{K}}^\varepsilon)^{*k}_{x,\mathrm{d}y}(n) = \mathcal{P}_\varepsilon(n\in\tau,\chi_1 > n, \varphi_{n-1}\in\mathrm{d}y|J_0=x).$$

Then, recalling that $q_{\varepsilon_c}(N+1) := \mathcal{P}_{\varepsilon_c}(\chi_1 = N+1)$, we can write, in analogy with (4.6),

$$\mathcal{P}_{\varepsilon_c}(\tau_{l_N} \le L, \tau_{r_N} \ge N-L, |J_{l_N}| \le L, |J_{r_N}| \le L | \chi_1 = N+1)$$
$$= \frac{1}{q_{\varepsilon_c}(N+1)} \cdot \sum_{a,b=0}^L \int_{x,y\in[-L,L],z\in\mathbb{R}} \widehat{\mathsf{U}}^{\varepsilon_c}_{0,\mathrm{d}x}(a) \cdot \widehat{\mathsf{K}}^{\varepsilon_c}_{x,\mathrm{d}y}(N+1-a-b)$$
$$\times \widehat{\mathsf{U}}^{\varepsilon_c}_{y,\mathrm{d}z}(b) \cdot \mathsf{K}^{\varepsilon_c}_{z,\{0\}}(1).$$

Applying relations (4.7) and (3.5), we obtain

(5.7)
$$\lim_{N\to\infty} \mathcal{P}_{\varepsilon_c}(\tau_{l_N} \le L, \tau_{r_N} \ge N-L, |J_{l_N}| \le L, |J_{r_N}| \le L | \chi_1 = N)$$
$$= \frac{\int_{x,y\in[-L,L],z\in\mathbb{R}} (\sum_a^L \widehat{\mathsf{U}}^{\varepsilon_c}_{0,\mathrm{d}x}(a))L^{\varepsilon_c}_{x,\mathrm{d}y}(\sum_{b=0}^L \widehat{\mathsf{U}}^{\varepsilon_c}_{y,\mathrm{d}z}(b))\mathsf{K}^{\varepsilon_c}_{z,\{0\}}(1)}{C_{\varepsilon_c}}.$$

However, the precise value of $C_{\varepsilon_c}$ is shown in [6], Section 7.3, to be

$$C_{\varepsilon_c} = ((1-\widehat{B}^{\varepsilon_c})^{-1} \circ L^{\varepsilon_c} \circ (1-\widehat{B}^{\varepsilon_c})^{-1} \circ \mathsf{K}^{\varepsilon_c})_{0,\{0\}},$$

where, of course, $\widehat{B}^\varepsilon_{x,\mathrm{d}y} := \sum_{n\in\mathbb{N}} \widehat{\mathsf{K}}^\varepsilon_{x,\mathrm{d}y}(n)$. Since $\sum_{n\in\mathbb{N}} \widehat{\mathsf{U}}^\varepsilon_{x,\mathrm{d}y}(n) = (1-\widehat{B}^\varepsilon)^{-1}_{x,\mathrm{d}y}$, by letting $L\to\infty$, the r.h.s. of (5.7) converges to 1 and equation (5.6) is proven.

**6. Proof of Theorem 1.6.** In this section we prove Theorem 1.6. We start discussing the topological and measurable structure of the space $\mathcal{M}([0,1])$ (for more details we refer to [12]).

6.1. *Finite signed measures.* We denote by $\mathcal{M}([0,1])$ the space of finite signed Borel measures on the interval $[0,1]$, that is of those set functions $\nu$ that can be written as $\nu = \nu_1 - \nu_2$, where $\nu_1$ and $\nu_2$ are finite nonnegative Borel measures on $[0,1]$ (since all the measures we deal with are Borel and finite, these adjectives will be dropped henceforth). According to the Hahn–Jordan decomposition [7], every $\nu \in \mathcal{M}([0,1])$ can be written in a unique way as $\nu = \nu^+ - \nu^-$, where $\nu^+$ and $\nu^-$ are nonnegative measures supported by *disjoint* Borel sets. Given $\nu \in \mathcal{M}([0,1])$, the nonnegative measure $|\nu| := \nu^+ + \nu^-$ is called the *total variation* of $\nu$. For $K \in \mathbb{R}^+$ we set

$$\mathcal{M}_K([0,1]) := \{\nu \in \mathcal{M}([0,1]) : |\nu|([0,1]) \le K\}.$$

Notice that $\mathcal{M}_K([0,1]) \subset \mathcal{M}_{K+1}([0,1])$ and that

$$\mathcal{M}([0,1]) = \bigcup_{K \in \mathbb{N}} \mathcal{M}_K([0,1]). \tag{6.1}$$

We recall that $C([0,1])$ denotes the space of continuous real functions defined on $[0,1]$. We equip the space $\mathcal{M}([0,1])$ with the topology of vague convergence, that, is the smallest topology on $\mathcal{M}([0,1])$ under which the map $\nu \mapsto \int f \, d\nu$ is continuous for every $f \in C([0,1])$, and with the corresponding Borel $\sigma$-field. We recall that $\nu_n \to \nu$ in $\mathcal{M}([0,1])$ if and only if $\int f \, d\nu_n \to \int f \, d\nu$ for all $f \in C([0,1])$ (see [9, 10] for a more explicit characterization).

The space $\mathcal{M}([0,1])$ is Hausdorff and separable [a dense countable subset is given by the measures $\sum_{i=1}^{n} a_i \delta_{b_i}(\cdot)$, for $n \in \mathbb{N}$ and $a_i, b_i \in \mathbb{Q}$]. The delicate point is that $\mathcal{M}([0,1])$ *is not metrizable*. However, we have the following result, proven in [12], Theorems 9.8.7 and 9.8.10.

LEMMA 6.1. *For every $K \in \mathbb{R}^+$, the space $\mathcal{M}_K([0,1])$ with the vague topology is compact and metrizable (and separable, hence Polish). Viceversa, if $A \subset \mathcal{M}([0,1])$ is compact, then $A \subset \mathcal{M}_K([0,1])$ for some $K \in \mathbb{N}$.*

By a *random signed measure* on $[0,1]$, we mean a random element $\boldsymbol{\nu}$, defined on some probability space $(\Omega, \mathcal{F}, P)$, and taking values in $\mathcal{M}([0,1])$. For instance, $\{\boldsymbol{\mu}_N\}_{N \in \mathbb{N}}$ defined in (1.18) (under the law $\mathbb{P}_{\varepsilon_c, N}$) is a sequence of random signed measures. For notational clarity, random signed measures will always be denoted by boldface symbols. The *law* of a random signed measure $\boldsymbol{\nu}$ is the probability measure $\boldsymbol{\nu} \circ P^{-1}$ on $\mathcal{M}([0,1])$. Given the random signed measures $\{\boldsymbol{\nu}_N\}_{N \in \mathbb{N}}$ and $\boldsymbol{\nu}$, we say that $\{\boldsymbol{\nu}_N\}_{N \in \mathbb{N}}$ converges in distribution on $\mathcal{M}([0,1])$ toward $\boldsymbol{\nu}$ if the law of $\boldsymbol{\nu}_N$ converges weakly to the law of $\boldsymbol{\nu}$, that is, for every bounded and continuous functional $F: \mathcal{M}([0,1]) \to \mathbb{R}$ we have $E[F(\boldsymbol{\nu}_N)] \to E[F(\boldsymbol{\nu})]$.

We are going to give sufficient conditions for convergence in distribution of random signed measures that will be applied in the next paragraphs. The path we follow is close to the standard one of proving tightness and checking the "convergence of the finite-dimensional distributions," but some additional care is required, due to the nonmetrizability of $\mathcal{M}([0,1])$. We recall that a sequence $\{\boldsymbol{\nu}_N\}_{N \in \mathbb{N}}$ of random signed measures on $[0,1]$ is said to be *tight* if for every $\delta > 0$ there exist a compact set $C \in \mathcal{M}([0,1])$ such that $P(\boldsymbol{\nu}_N \in C) \geq 1 - \delta$ for large $N$. Equivalently, $\{\boldsymbol{\nu}_N\}_{N \in \mathbb{N}}$ is tight if and only if for every $\delta > 0$ there exist $K, N_0 \in \mathbb{N}$ such that

$$P(|\boldsymbol{\nu}_N|([0,1]) \leq K) \geq 1 - \delta \qquad \forall N \geq N_0. \tag{6.2}$$

Although $\mathcal{M}([0,1])$ is not Polish, the first half of Prohorov's Theorem still holds:

LEMMA 6.2. *If the sequence of random signed measures $\{\boldsymbol{\nu}_N\}_{N\in\mathbb{N}}$ is tight, then there is a subsequence $\{\boldsymbol{\nu}_{N_k}\}_{k\in\mathbb{N}}$ which converges in distribution on $\mathcal{M}([0,1])$.*

The proof of this lemma is given in Appendix B.2. Next, for $t \in [0,1]$ we define the measurable map $F_t : \mathcal{M}([0,1]) \to \mathbb{R}$ by

$$F_t(\nu) := \nu([0,t]).$$

For $k \in \mathbb{N}$ and $0 \le a_1 \le \cdots \le a_k \le 1$, a random signed measure $\boldsymbol{\nu}$ determines the law on $\mathbb{R}^k$ defined by

$$(F_{a_1}(\boldsymbol{\nu}),\ldots,F_{a_k}(\boldsymbol{\nu})) \circ P^{-1} = (\boldsymbol{\nu}([0,a_1]),\boldsymbol{\nu}([0,a_2]),\ldots,\boldsymbol{\nu}([0,a_k])) \circ P^{-1},$$

where $P$ is the underlying probability measure. These laws are called the *finite dimensional distributions* of the random signed measure $\boldsymbol{\nu}$. Notice that if $\boldsymbol{\nu}_1$ and $\boldsymbol{\nu}_2$ have the same finite dimensional distributions then they have the same law on $\mathcal{M}([0,1])$, because the $\sigma$-field generated by the maps $\{F_t\}_{t\in[0,1]}$ coincides with the Borel $\sigma$-field of $\mathcal{M}([0,1])$. In other terms, the finite dimensional distributions determine laws on $\mathcal{M}([0,1])$.

We are ready to put together tightness and convergence of the finite-dimensional distributions to yield convergence in distribution on $\mathcal{M}([0,1])$. The next proposition, proven in Appendix B.1, is sufficient for our purposes.

PROPOSITION 6.3. *Let $\{\boldsymbol{\nu}_N\}_{N\in\mathbb{N}}$ be a* tight *sequence of random signed measures on $[0,1]$. Assume that the finite-dimensional distributions of $\boldsymbol{\nu}_N$ converge, that is, $\forall k \in \mathbb{N}$ and for all $0 < a_1 < \cdots < a_k < 1$ there is a probability measure $\lambda^{(k)}_{a_1,\ldots,a_k}(\cdot)$ on $\mathbb{R}^k$ such that*

(6.3) $$(\boldsymbol{\nu}_N([0,a_1]),\boldsymbol{\nu}_N([0,a_2]),\ldots,\boldsymbol{\nu}_N([0,a_k])) \xrightarrow{d} \lambda^{(k)}_{a_1,\ldots,a_k} \qquad (N \to \infty).$$

*Assume, moreover, that for every $x \in [0,1]$ and $\eta > 0$*

(6.4) $$\lim_{\delta \to 0} \limsup_{N\to\infty} P(|\boldsymbol{\nu}_N|([x-\delta,x+\delta]) > \eta) = 0.$$

*Then $\{\boldsymbol{\nu}_N\}_{N\in\mathbb{N}}$ converges in distribution on $\mathcal{M}([0,1])$ toward a random signed measure whose finite-dimensional distributions are $\lambda^{(k)}_{a_1,\ldots,a_k}$.*

The reason for requiring the extra condition (6.4) is that the map $F_t$ is not continuous on $\mathcal{M}([0,1])$ and, therefore, the convergence in distribution on $\mathcal{M}([0,1])$ does not imply automatically the convergence of the finite-dimensional distributions.

6.2. *Preparation.* Remember the definition (1.18) of the random signed measure $\boldsymbol{\mu}_N$ under $\mathbb{P}_{\varepsilon_c,N}$ that we look at as a random element of the space $\mathcal{M}([0,1])$. Our goal is to show that $\boldsymbol{\mu}_N$ under $\mathbb{P}_{\varepsilon_c,N}$ converges in distribution as $N \to \infty$ toward the random measure $dL$, defined in Section 1.4, using Proposition 6.3.

We start restating for $\boldsymbol{\mu}_N$ the convergence of the finite-dimensional distributions and the extra-condition (6.4), which are interesting by themselves.

THEOREM 6.4. *For every $k \in \mathbb{N}$ and for all $a_1, \ldots, a_k \in (0,1)$ with $a_i \leq a_{i+1}$, $i = 1, \ldots, k$, we have as $N \to \infty$*

$$(6.5) \quad (\boldsymbol{\mu}_N((0, a_1]), \boldsymbol{\mu}_N((a_1, a_2]), \ldots, \boldsymbol{\mu}_N((a_{k-1}, a_k]))$$
$$under \quad \mathbb{P}_{\varepsilon_c, N} \xrightarrow{d} (L_{a_1}, L_{a_2} - L_{a_1}, \ldots, L_{a_k} - L_{a_{k-1}}),$$

*where $\xrightarrow{d}$ denotes convergence in distribution on $\mathbb{R}^k$. Moreover, $\forall x \in [0,1]$, $\forall \eta > 0$,*

$$(6.6) \quad \lim_{\delta \to 0} \limsup_{N \to \infty} \mathbb{P}_{\varepsilon_c, N}(|\boldsymbol{\mu}_N|([x - \delta, x + \delta]) > \eta) = 0.$$

Notice that the vectors in (6.5) differ from those in (6.3) just by a linear transformation, because it is simpler to work with $\boldsymbol{\mu}_N((a_{i-1}, a_i])$ than with $\boldsymbol{\mu}_N((0, a_i]) = \boldsymbol{\mu}_N([0, a_i])$.

The proof of Theorem 6.4 is given in Section 6.3, while the tightness of the sequence $\{\boldsymbol{\mu}_N\}_N$ under $\mathbb{P}_{\varepsilon_c,N}$ is proven in Section 6.4. Thanks to Proposition 6.3, this completes the proof of Theorem 1.6. The rest of this paragraph is devoted to a basic lemma.

LEMMA 6.5. *Fix any $\delta \in (0,1)$. Given any sequence of events $\{B_N\}_{N \in \mathbb{N}}$ such that $B_N \in \sigma(\{\varphi_i\}_{0 \leq i \leq \delta N})$, that is, $B_N$ depends on the field of length $\delta N$, the following relation holds:*

$$\mathbb{P}_{\varepsilon_c, N}(B_N) = \mathcal{P}_{\varepsilon_c}(B_N) + o(1) \qquad (N \to \infty).$$

PROOF. Thanks to relation (2.18), it suffices to prove that

$$(6.7) \quad \mathcal{P}_{\varepsilon_c}(B_N | N + 1 \in \chi) = \mathcal{P}_{\varepsilon_c}(B_N) + o(1) \qquad (N \to \infty).$$

Introducing the variable $\xi_\delta := \min\{\chi \cap [\delta N, \infty)\} - \max\{\chi \cap [0, \delta N]\}$, we claim that

$$(6.8) \quad \mathcal{P}_{\varepsilon_c}\left(\xi_\delta \geq \frac{N}{\log N} \Big| N + 1 \in \chi\right) = o(1), \qquad \mathcal{P}_{\varepsilon_c}\left(\xi_\delta \geq \frac{N}{\log N}\right) = o(1).$$

In fact, these relations are proven in Appendix A with explicit bounds [cf. (A.7)–(A.9) and (A.12)] in the special case $\delta = \frac{1}{2}$, but the proof carries over to the general case with no change. We introduce the variable $d_\delta := \min\{\chi \cap [\delta N, \infty)\} - \lfloor \delta N \rfloor$, and we note that $d_\delta \leq \xi_\delta$. Thanks to (6.8), we can rephrase (6.7) as

$$\mathcal{P}_{\varepsilon_c}\left(B_N, d_\delta \leq \frac{N}{\log N} \bigg| N + 1 \in \chi\right) = \mathcal{P}_{\varepsilon_c}\left(B_N, d_\delta \leq \frac{N}{\log N}\right) + o(1)$$
(6.9)
$$(N \to \infty).$$

We recall from Section 3.1 that the process $\{\chi_n\}_n$ under $\mathcal{P}_{\varepsilon_c}$ is a renewal process with step law $q_{\varepsilon_c}(n) = \mathcal{P}_{\varepsilon_c}(\chi_1 = n)$. Denoting by $u_{\varepsilon_c}(n) := \sum_{k \geq 0} q_{\varepsilon_c}^{*k}(n)$ the corresponding renewal mass function, we can write the l.h.s. of (6.9) as

$$\mathcal{P}_{\varepsilon_c}\left(B_N, d_\delta \leq \frac{N}{\log N} \bigg| N + 1 \in \chi\right)$$
$$= \sum_{k=0}^{\lfloor N/\log N \rfloor} \mathcal{P}_{\varepsilon_c}(B_N, d_\delta = k) \cdot \frac{u_{\varepsilon_c}(N + 1 - \lfloor \delta N \rfloor - k)}{u_{\varepsilon_c}(N + 1)}.$$

Since $q_{\varepsilon_c}(n) \sim C_{\varepsilon_c}/n^2$ as $n \to \infty$ [see (3.5)], by Theorem 8.7.5 of [3], we have $u_{\varepsilon_c}(n) \sim 1/(C_{\varepsilon_c} \log n)$. Therefore, $u_{\varepsilon_c}(N + 1 - \lfloor \delta N \rfloor - k)/u_{\varepsilon_c}(N + 1) = 1 + o(1)$ as $N \to \infty$, uniformly for $k$ in the range of summation, and (6.9) is proven. □

COROLLARY 6.6. *To prove equations (6.5) and (6.6), one can replace the law* $\mathbb{P}_{\varepsilon_c, N}$ *by* $\mathcal{P}_{\varepsilon_c}$.

6.3. *Proof of Theorem 6.4.* We introduce the sequences $\{A_k\}_{k \in \mathbb{N}}$ and $\{\widetilde{A}_k\}_{k \in \mathbb{N}}$ that give the area respectively under the processes $\{\varphi_i\}_i$ and $\{|\varphi_i|\}_i$ between two consecutive adjacent contact points:

$$(6.10) \qquad A_k := \sum_{i=\chi_{k-1}+1}^{\chi_k} \varphi_i, \qquad \widetilde{A}_k := \sum_{i=\chi_{k-1}+1}^{\chi_k} |\varphi_i|.$$

We also introduce the corresponding partial sum processes:

$$(6.11) \qquad S_n := A_1 + \cdots + A_n, \qquad \widetilde{S}_n := \widetilde{A}_1 + \cdots + \widetilde{A}_n.$$

Note that the variables $\{A_k\}_{k \in \mathbb{N}}$ are i.i.d. under $\mathcal{P}_{\varepsilon_c}$ and, hence, $\{S_n\}_{n \geq 0}$ is a real random walk, and analogous statements hold for $\{\widetilde{A}_k\}_{k \in \mathbb{N}}$ and $\{\widetilde{S}_n\}_{n \geq 0}$. In fact, the epochs $\{\chi_k\}_{k \geq 0}$ cut the field into independent segments, because $\{\chi_k\}_{k \geq 0}$ under $\mathcal{P}_{\varepsilon_c}$ is a genuine renewal process [cf. Section 3.1] and, furthermore, the excursions $\{e_k\}_{k \in \mathbb{N}}$ are independent conditionally on $\{(\tau_k, J_k)\}_{k \in \mathbb{Z}^+}$; compare Section 2.3.

The crucial fact is that the random walk $\{S_n\}_n$ under $\mathcal{P}_{\varepsilon_c}$ is in the domain of attraction of the symmetric stable Lévy process of index $2/5$ and, analogously, $\{\widetilde{S}_n\}_n$ is in the domain of attraction of the stable subordinator of index $2/5$. In fact, we have the following:

PROPOSITION 6.7.  *There exist positive constants $\mathcal{C}, \widetilde{\mathcal{C}}$ such that*

$$(6.12) \qquad \mathcal{P}_{\varepsilon_c}(A_1 > x) \sim \frac{\mathcal{C}}{x^{2/5}}, \qquad \mathcal{P}_{\varepsilon_c}(\widetilde{A}_1 > x) \sim \frac{\widetilde{\mathcal{C}}}{x^{2/5}} \qquad (x \to +\infty).$$

PROOF.  By Lemma 5.1 and the Continuous Mapping Theorem, as $n \to \infty$, we have that

$$(6.13) \qquad \int_0^1 \widehat{\varphi}_n(t)\,dt = \frac{1}{\sigma}\frac{1}{n^{5/2}} \sum_{i=1}^n \varphi_i \quad \text{under} \quad \mathcal{P}_{\varepsilon_c}(\cdot|\chi_1 = n) \xrightarrow{d} \int_0^1 \widehat{I}_t\,dt,$$

where $\xrightarrow{d}$ denotes convergence in distribution on $\mathbb{R}$ and the process $\{\widehat{I}_t\}_{t \in [0,1]}$ was introduced in (1.11). Note that $\int_0^1 \widehat{I}_t\,dt$ is a Gaussian random variable, whose variance equals $\frac{1}{720}$ (see Appendix B.3), hence,

$$(6.14) \qquad \Phi(z) := P\left(\int_0^1 \widehat{I}_t\,dt > z\right) = \frac{6\sqrt{10}}{\sqrt{\pi}} \int_z^\infty e^{-360 t^2}\,dt.$$

For $z \in \mathbb{R}$ and $n \in \mathbb{N}$, we set

$$(6.15) \qquad \Phi_n(z) := \mathcal{P}_{\varepsilon_c}\left(\frac{A_1}{\sigma n^{5/2}} > z \,\Big|\, \chi_1 = n\right),$$

and note that equation (6.13) yields $\Phi_n(z) \to \Phi(z)$ as $n \to \infty$, for every $z \in \mathbb{R}$.

Recalling the notation $q_{\varepsilon_c}(n) = \mathcal{P}_{\varepsilon_c}(\chi_1 = n)$, we can write

$$\mathcal{P}_{\varepsilon_c}(A_1 > x) = \sum_{n \in \mathbb{N}} q_{\varepsilon_c}(n) \mathcal{P}_{\varepsilon_c}(A_1 > x | \chi_1 = n)$$

$$= \sum_{n \in \mathbb{N}} q_{\varepsilon_c}(n) \Phi_n\left(\frac{x}{\sigma n^{5/2}}\right).$$

Let us rewrite the r.h.s. above by putting in evidence the factor $s := \frac{n\sigma^{2/5}}{x^{2/5}}$:

$$(6.16) \qquad \mathcal{P}_{\varepsilon_c}(A_1 > x) = \frac{1}{x^{2/5}}\left\{\sigma^{2/5} \cdot \frac{\sigma^{2/5}}{x^{2/5}} \sum_{s \in \sigma^{2/5}/n^{2/5}\mathbb{N}} \left[\frac{x^{4/5}}{\sigma^{4/5}} q\left(\frac{x^{2/5}}{\sigma^{2/5}} s\right)\right] \right.$$
$$\left. \times \Phi_{x^{2/5}/\sigma^{2/5} s}\left(\frac{1}{s^{5/2}}\right)\right\}.$$

Since $\Phi_n(z) \to \Phi(z)$ and $q_{\varepsilon_c}(n) \sim C_{\varepsilon_c}/n^2$ as $n \to \infty$ [cf. (3.5)], for every $s > 0$ we have

$$\Phi_{x^{2/5}/\sigma^{2/5}s}\left(\frac{1}{s^{5/2}}\right) \longrightarrow \Phi\left(\frac{1}{s^{5/2}}\right),$$

$$\frac{x^{4/5}}{\sigma^{4/5}} q\left(\frac{x^{2/5}}{\sigma^{2/5}} s\right) \longrightarrow \frac{C_{\varepsilon_c}}{s^2} \qquad (x \to +\infty).$$

Moreover, we claim that the following bound holds true (see below):

$$(6.17) \qquad \mathcal{P}_{\varepsilon_c}(A_1 > x | \chi_1 = n) \leq (const.) \frac{n^5}{x^2}.$$

Then a Riemann-sum argument shows that the term in brackets in (6.16) does converge toward the corresponding integral, that is, as $x \to \infty$,

$$\sigma^{2/5} \cdot \frac{\sigma^{2/5}}{x^{2/5}} \sum_{s \in \sigma^{2/5}/n^{2/5}\mathbb{N}} \left[\frac{x^{4/5}}{\sigma^{4/5}} q\left(\frac{x^{2/5}}{\sigma^{2/5}} s\right)\right] \cdot \Phi_{x^{2/5}/\sigma^{2/5}s}\left(\frac{1}{s^{5/2}}\right)$$

$$(6.18) \qquad \longrightarrow \sigma^{2/5} \int_0^\infty \frac{C_{\varepsilon_c}}{s^2} \Phi\left(\frac{1}{s^{5/2}}\right) ds$$

$$= C_{\varepsilon_c} \sigma^{2/5} \frac{6\sqrt{10}}{\sqrt{\pi}} \int_0^\infty t^{2/5} e^{-360 t^2} \, dt =: \mathcal{C},$$

having used (6.14). This proves the first relation in (6.12), with an explicit formula for $\mathcal{C}$.

The variable $\widetilde{A}_1$ is treated in a similar way. In fact, in analogy with (6.13), Lemma 5.1 and the Continuous Mapping Theorem yield, as $n \to \infty$,

$$(6.19) \qquad \frac{\widetilde{A}_1}{\sigma n^{5/2}} \quad \text{under} \quad \mathcal{P}_{\varepsilon_c}(\cdot | \chi_1 = n) \xrightarrow{d} \int_0^1 |\widehat{I}_t| \, dt,$$

and, moreover, the following bound holds (see below):

$$(6.20) \qquad \mathcal{P}_{\varepsilon_c}(\widetilde{A}_1 > x | \chi_1 = n) \leq (const.) \frac{n^5}{x^2}.$$

Then, arguing exactly as above, a Riemann-sum approximation shows that the second relation in (6.12) holds true, with

$$(6.21) \quad \widetilde{\mathcal{C}} := \sigma^{2/5} \int_0^\infty \frac{C_{\varepsilon_c}}{s^2} \widetilde{\Phi}\left(\frac{1}{s^{5/2}}\right) ds = \sigma^{2/5} \int_0^\infty \frac{2}{5} \frac{C_{\varepsilon_c}}{t^{3/5}} \widetilde{\Phi}(t) \, dt < \infty,$$

where, of course, $\widetilde{\Phi}(t) := P(\int_0^1 |\widehat{I}_s| \, ds > t)$.

To complete the proof, it remains to prove (6.20), which implies (6.17), because $A_1 \leq \widetilde{A}_1$. To this purpose, we exploit the Brascamp–Lieb inequality. We recall from Section 3.3 that the law of the vector $(Z_1, \ldots, Z_n)$ under $\mathbf{P}^{(0,0)}$ has the form $\mu_H(dx) = e^{-H(x)} \, dx$, $x \in \mathbb{R}^n$, where $H(x) = \frac{1}{2} A(x) + R(x)$ and $A(\cdot), R(\cdot)$ are defined in (3.10). Fixing $m \leq n$ and $t_1, \ldots, t_m \in \{1, \ldots, n\}$, the law $\mu_H^*(\mathrm{d}x) := \mu_H(\mathrm{d}x | x_{t_1} = 0, \ldots, x_{t_m} = 0)$ satisfies the Brascamp–Lieb inequality (3.9): choosing $\Gamma(x) = x_k^2$, with $1 \leq k \leq n$, we obtain

$$\mathbf{E}^{(0,0)}(Z_k^2 | Z_{t_1} = 0, \ldots, Z_{t_m} = 0) = E_H^*(x_k^2) \leq E_A(x_k^2) = \frac{k(k+1)(2k+1)}{6\gamma},$$

where we observe that $E_H^*(x_k) = 0$ by symmetry and the last equality is just a straightforward Gaussian computation, because $\mu_A$ is nothing but the law of the integral of a random walk with Gaussian steps $\sim \mathcal{N}(0, \gamma^{-1})$ [cf. (2.2) and (2.4)]. Setting $\mathbf{P}_*^{(0,0)}(\cdot) := \mathbf{P}^{(0,0)}(\cdot | Z_{t_1} = 0, \ldots, Z_{t_m} = 0)$ for conciseness and using the Chebyshev and Cauchy–Schwarz inequalities, we obtain

$$\mathbf{P}^{(0,0)}\left(\sum_{k=1}^n |Z_k| > x \Big| Z_{t_1} = 0, \ldots, Z_{t_m} = 0\right)$$

(6.22)
$$\leq \frac{1}{x^2} \sum_{k,l=1}^n \mathbf{E}_*^{(0,0)}(|Z_k| \cdot |Z_l|) \leq \frac{1}{x^2}\left(\sum_{k=1}^n (\mathbf{E}_*^{(0,0)}(Z_k^2))^{1/2}\right)^2$$

$$\leq \frac{1}{\gamma x^2}\left(\sum_{k=1}^n k^{3/2}\right)^2 \leq (const.) \frac{n^5}{x^2}.$$

Now observe that $\mathbb{P}_{\varepsilon_c, n-1}(\cdot) = \mathcal{P}_{\varepsilon_c}(\cdot | n \in \chi)$ [cf. (2.18)], hence, $\mathbb{P}_{\varepsilon_c, n-1}(\cdot | \chi_1 = n) = \mathcal{P}_{\varepsilon_c}(\cdot | \chi_1 = n)$. We set

$$\mathcal{A}_{n-2} := \{A \subset \{1, n-2\} : 1 \notin A, n-2 \notin A, \{l, l+1\} \not\subset A, \ \forall 1 \leq l \leq n-3\},$$

and we use the notation $\tau_{[1, n-2]} := \tau \cap [1, n-2]$. Noting that $\mathcal{A}_{n-2}$ represents all the possible values of the variable $\tau_{[1, n-2]}$ under $\mathcal{P}_{\varepsilon_c}(\cdot | \chi_1 = n)$, we can write

$$\mathcal{P}_{\varepsilon_c}\left(\sum_{k=1}^n |\varphi_k| > x \Big| \chi_1 = n\right) = \mathbb{P}_{\varepsilon_c, n-1}\left(\sum_{k=1}^n |\varphi_k| > x \Big| \chi_1 = n\right)$$

$$= \sum_{A \in \mathcal{A}_{n-2}} \mathbb{P}_{\varepsilon_c, n-1}\left(\sum_{k=1}^n |\varphi_k| > x \Big| \tau_{[1, n-2]} = A\right)$$

$$\times \mathbb{P}_{\varepsilon_c, n-1}(\tau_{[1, n-2]} = A | \chi_1 = n).$$

However, combining (2.5) with (6.22), we have

$$\mathbb{P}_{\varepsilon_c, n-1}\left(\sum_{k=1}^n |\varphi_k| > x \Big| \tau_{[1, n-2]} = A\right)$$

$$= \mathbf{P}^{(0,0)}\left(\sum_{k=1}^n |Z_k| > x \Big| Z_i = 0, \ \forall i \in A \cup \{n-1, n\}\right) \leq (const.) \frac{n^5}{x^2},$$

and the proof of (6.20) is completed. □

We now denote by $\{\widetilde{L}_t\}_{t\in[0,1]}$ the stable subordinator of index $\frac{2}{5}$, normalized so that its Lévy measure equals $\frac{2}{5}\widetilde{\mathcal{C}}x^{-2/5-1}\mathbf{1}_{\{x>0\}}\,dx$, so that $P(\widetilde{L}_1 > x) \sim \widetilde{\mathcal{C}}x^{-2/5}$ as $x \to \infty$. By Proposition 6.7, we have $\mathcal{P}_{\varepsilon_c}(\widetilde{A}_1 > x) \sim P(\widetilde{L}_1 > x)$ as $x \to \infty$, hence, by the standard theory of stability ([8], Chapter XVII.5), $\widetilde{A}_1$ is in the domain of attraction of $\widetilde{L}_1$ and we have

$$(6.23) \quad \frac{1}{n^{5/2}}\widetilde{S}_n = \frac{1}{n^{5/2}}\sum_{i=1}^n \widetilde{A}_i \quad \text{under} \quad \mathcal{P}_{\varepsilon_c} \xrightarrow{d} \widetilde{L}_1 \qquad (n \to \infty).$$

Next let $\{L_t\}_{t\in[0,1]}$ be the symmetric stable Lévy process of index $\frac{2}{5}$, with Lévy measure given by $c_L|x|^{-2/5-1}dx$, where $c_L := \mathcal{C}/C_{\varepsilon_c}$ [we recall that $C_{\varepsilon_c}$ is the constant appearing in (3.5)]. In particular, we have $P(L_1 > x) = P(L_1 < -x) \sim c_L x^{-2/5}$ as $x \to \infty$. Then Proposition 6.7 yields $\mathcal{P}_{\varepsilon_c}(A_1 > x) \sim P((C_{\varepsilon_c})^{5/2}L_1 > x)$ as $x \to \infty$, and since $\mathcal{P}_{\varepsilon_c}(A_1 > x) = \mathcal{P}_{\varepsilon_c}(A_1 < -x)$ by symmetry, it follows by the theory of stability that

$$(6.24) \quad \frac{1}{n^{5/2}}S_n = \frac{1}{n^{5/2}}\sum_{i=1}^n A_i \quad \text{under} \quad \mathcal{P}_{\varepsilon_c} \xrightarrow{d} (C_{\varepsilon_c})^{5/2}L_1 \qquad (n \to \infty).$$

Notice that by (6.18) the constant $c_L := \mathcal{C}/C_{\varepsilon_c}$ equals

$$(6.25) \quad c_L = \frac{6\sqrt{10}}{\sqrt{\pi}}\sigma^{2/5}\int_0^\infty s^{2/5}e^{-360s^2}\,ds = \frac{3\sqrt{10}}{\sqrt{\pi}(360)^{7/10}}\Gamma\left(\frac{7}{10}\right)\sigma^{2/5},$$

where $\Gamma(x) := \int_0^\infty t^{x-1}e^{-t}\,dt$ is the usual Gamma function and the second equality follows by a simple change of variables. We also recall that $L_t \stackrel{d}{=} t^{5/2}L_1$ and $\widetilde{L}_t \stackrel{d}{=} t^{5/2}\widetilde{L}_1$.

We are ready to prove (6.6), with the law $\mathbb{P}_{\varepsilon_c,N}$ replaced by $\mathcal{P}_{\varepsilon_c}$, thanks to Corollary 6.6. It is convenient to extend the definition of $\iota_N$ to a noninteger argument, by setting $\iota[t] := \sup\{k \in \mathbb{Z}^+ : \chi_k \leq t\}$ for $t \in \mathbb{R}$; compare (3.3). By the definitions (1.17) and (1.18) of $\widetilde{\varphi}_N$ and $\boldsymbol{\mu}_N$, we immediately obtain the following upper bound:

$$(6.26) \quad |\boldsymbol{\mu}_N|([x-\delta, x+\delta]) \leq \left(\frac{\log N}{N}\right)^{5/2}(\widetilde{S}_{\iota[(x+\delta)N]+1} - \widetilde{S}_{\iota[(x-\delta)N]}).$$

Since $\{\chi_k\}_{k\geq 0}$ is a genuine renewal process with $\mathcal{P}_{\varepsilon_c}(\chi_1 = n) \sim C_{\varepsilon_c}/n^2$, Theorem 8.8.1 of [3] yields $\chi_k/(k\log k) \to C_{\varepsilon_c}$ as $k \to \infty$, $\mathcal{P}_{\varepsilon_c}$-a.s., and since $\chi_{\iota[t]} \leq t \leq \chi_{\iota[t]+1}$, it follows that

$$(6.27) \quad \frac{\iota[t]}{t/\log t} \longrightarrow \frac{1}{C_{\varepsilon_c}} \qquad \text{as } t \to \infty, \mathcal{P}_{\varepsilon_c}\text{-a.s.}$$

Therefore, for every $\kappa > 0$ we can choose $N$ sufficiently large such that

$$\mathcal{P}_{\varepsilon_c}\left(\left\{\iota[(x+\delta)N]+1 > \frac{1}{C_{\varepsilon_c}}\frac{(x+2\delta)N}{\log N}\right\} \cup \left\{\iota[(x-\delta)N] < \frac{1}{C_{\varepsilon_c}}\frac{(x-2\delta)N}{\log N}\right\}\right)$$
$$\leq \kappa.$$

Then by (6.26) for any $\eta > 0$ and for large $N$, we can write

$$\mathcal{P}_{\varepsilon_c}(|\boldsymbol{\mu}_N|([x-\delta, x+\delta]) > \eta)$$
$$\leq \kappa + \mathcal{P}_{\varepsilon_c}\left(\left(\frac{\log N}{N}\right)^{5/2}(\widetilde{S}_{\lfloor (x+2\delta)N/(C_{\varepsilon_c}\log N)\rfloor} - \widetilde{S}_{\lfloor (x-2\delta)N/(C_{\varepsilon_c}\log N)\rfloor}) > \eta\right).$$

However, for $a, b \in \mathbb{N}$ with $a \leq b$ we have $\widetilde{S}_b - \widetilde{S}_a \stackrel{d}{=} \widetilde{S}_{b-a}$. Then letting $N \to \infty$ and recalling (6.23), we have

$$\limsup_{N\to\infty} \mathcal{P}_{\varepsilon_c}(|\boldsymbol{\mu}_N|([x-\delta, x+\delta]) > \eta) \leq \kappa + P\left(\left(\frac{4\delta}{C_{\varepsilon_c}}\right)^{5/2}\widetilde{L}_1 > \eta\right).$$

Letting $\delta \to 0$, the last term vanishes and since $\kappa$ was arbitrary, equation (6.6) is proven.

Next we prove (6.5), again with the law $\mathbb{P}_{\varepsilon_c,N}$ replaced by $\mathcal{P}_{\varepsilon_c}$, thanks to Corollary 6.6. We claim that (6.5) is equivalent to the following relation:

$$\left(\frac{\log N}{N}\right)^{5/2}(S_{\lfloor a_1 N/(C_{\varepsilon_c}\log N)\rfloor}, (S_{\lfloor a_2 N/(C_{\varepsilon_c}\log N)\rfloor} - S_{\lfloor a_1 N/(C_{\varepsilon_c}\log N)\rfloor}),$$
(6.28)
$$\ldots, (S_{\lfloor a_k N/(C_{\varepsilon_c}\log N)\rfloor} - S_{\lfloor a_{k-1}N/(C_{\varepsilon_c}\log N)\rfloor}))$$
$$\text{under } \mathcal{P}_{\varepsilon_c} \xrightarrow{d} (L_{a_1}, L_{a_2} - L_{a_1}, \ldots, L_{a_k} - L_{a_{k-1}}).$$

To prove the claim, it suffices to show that the difference between the vectors in the first lines of (6.5) and (6.28) converges in $\mathcal{P}_{\varepsilon_c}$-probability to zero as $N \to \infty$. It is sufficient to focus on each component: so we need to prove that

(6.29)
$$\lim_{N\to\infty} \mathcal{P}_{\varepsilon_c}\left(\left|\boldsymbol{\mu}_N((a,b]) - \left(\frac{\log N}{N}\right)^{5/2}(S_{\lfloor bN/(C_{\varepsilon_c}\log N)\rfloor} - S_{\lfloor aN/(C_{\varepsilon_c}\log N)\rfloor})\right| \geq \eta\right)$$
$$= 0,$$

for every $\eta > 0$ and for all $a, b \in [0, 1)$ with $a < b$. Fix $\delta > 0$ and observe that, by (6.27),

$$\lim_{N\to\infty} \mathcal{P}_{\varepsilon_c}\left(\iota[aN] \in \left\lfloor\frac{aN}{C_{\varepsilon_c}\log N}\right\rfloor \cdot (1-\delta, 1+\delta),\right.$$
$$\left.\iota[bN] \in \left\lfloor\frac{bN}{C_{\varepsilon_c}\log N}\right\rfloor \cdot (1-\delta, 1+\delta)\right) = 1.$$

Therefore, we can restrict ourselves on this event, where, using the definitions (1.17) and (1.18) of $\widetilde{\varphi}_N$ and $\boldsymbol{\mu}_N$, we can write

$$\left| \boldsymbol{\mu}_N((a,b]) - \left(\frac{\log N}{N}\right)^{5/2} (S_{\lfloor bN/(C_{\varepsilon_c} \log N) \rfloor} - S_{\lfloor aN/(C_{\varepsilon_c} \log N) \rfloor}) \right|$$

$$\leq \left(\frac{\log N}{N}\right)^{5/2} \{(\widetilde{S}_{\lfloor a(1+\delta)N/(C_{\varepsilon_c} \log N) \rfloor} - \widetilde{S}_{\lfloor a(1-\delta)N/(C_{\varepsilon_c} \log N) \rfloor})$$

$$+ (\widetilde{S}_{\lfloor b(1+\delta)N/(C_{\varepsilon_c} \log N) \rfloor} - \widetilde{S}_{\lfloor b(1-\delta)N/(C_{\varepsilon_c} \log N) \rfloor})\}.$$

However, $(\widetilde{S}_b - \widetilde{S}_a) + (\widetilde{S}_d - \widetilde{S}_c) \stackrel{d}{=} \widetilde{S}_{(b-a)+(d-c)}$ for $a \leq b \leq c \leq d$, and as $N \to \infty$ by (6.23), we have

$$\mathcal{P}_{\varepsilon_c}\left(\left(\frac{\log N}{N}\right)^{5/2} \widetilde{S}_{\lfloor (a+b) \cdot 2\delta N/(C_{\varepsilon_c} \log N) \rfloor} \geq \eta\right) \longrightarrow P\left(\frac{(a+b) \cdot 2\delta}{C_{\varepsilon_c}} L_1 \geq \eta\right).$$

The last term vanishes as $\delta \to 0$, hence, (6.29) is proven.

It finally remains to prove equation (6.28). Both the vector in the l.h.s. and the one in the r.h.s. of that equation have independent components, therefore, it suffices to prove the convergence of each component, that is, that for every $a \in (0,1)$ as $N \to \infty$,

$$\left(\frac{\log N}{N}\right)^{5/2} S_{\lfloor aN/(C_{\varepsilon_c} \log N) \rfloor} = \left(\frac{\log N}{N}\right)^{5/2} \sum_{i=1}^{\lfloor aN/(C_{\varepsilon_c} \log N) \rfloor} A_i$$

$$\text{under } \mathcal{P}_{\varepsilon_c} \xrightarrow{d} L_a.$$

However, recalling that $L_a \stackrel{d}{=} a^{5/2} L_1$, this relation follows immediately from (6.24), so that the proof of Theorem 6.4 is completed.

6.4. *Tightness of $\{\mu_N\}_N$.* We finally prove the tightness of the sequence $\{\boldsymbol{\mu}_N\}_{N \in \mathbb{N}}$, that is, for every $\delta > 0$ there exist $K, N_0 \in \mathbb{N}$ such that

(6.30) $\quad \mathbb{P}_{\varepsilon_c,N}(|\boldsymbol{\mu}_N|([0,1]) \leq K) \geq 1 - \delta \qquad \forall N \geq N_0.$

Since $\boldsymbol{\mu}_N(\{\frac{1}{2}\}) = 0$, we can write $\boldsymbol{\mu}_N([0,1]) = \boldsymbol{\mu}_N([0,\frac{1}{2}]) + \boldsymbol{\mu}_N([\frac{1}{2},1])$. However, by symmetry, $\boldsymbol{\mu}_N([0,\frac{1}{2}]) \stackrel{d}{=} \boldsymbol{\mu}_N([\frac{1}{2},1])$ under $\mathbb{P}_{\varepsilon_c,N}$, hence, it suffices to show that

$$\mathbb{P}_{\varepsilon_c,N}(|\boldsymbol{\mu}_N|([0,\tfrac{1}{2}]) \leq K/2) \geq 1 - \delta/2 \qquad \forall N \geq N_0.$$

Now notice that the event $\{|\boldsymbol{\mu}_N|([0,\frac{1}{2}]) \leq \frac{K}{2}\}$ belongs to the $\sigma$-field $\sigma(\{\varphi_i\}_{0 \leq i \leq N/2})$, hence, we can apply Lemma 6.5 and we are left with showing that for every $\delta > 0$ there exist $K, N_0 \in \mathbb{N}$ such that

(6.31) $\quad \mathcal{P}_{\varepsilon_c}(|\boldsymbol{\mu}_N|([0,\tfrac{1}{2}]) \leq K/2) \geq 1 - \delta/4 \qquad \forall N \geq N_0.$

We recall that $\iota[t] := \sup\{k \in \mathbb{Z}^+ : \chi_k \leq t\}$, for $t \in \mathbb{R}$. From the definitions (1.17), (1.18) and (6.11) of $\widetilde{\varphi}_N$, $\boldsymbol{\mu}_N$ and $\widetilde{S}_n$ respectively, the inclusion bound yields

$$\mathcal{P}_{\varepsilon_c}\left(|\boldsymbol{\mu}_N|\left(\left[0, \frac{1}{2}\right]\right) > K/2\right) \leq \mathcal{P}_{\varepsilon_c}\left(\left(\frac{\log N}{N}\right)^{5/2} \widetilde{S}_{\iota[N/2]+1} > \frac{K}{2}\right)$$

$$\leq \mathcal{P}_{\varepsilon_c}\left(\iota[N/2] + 1 > \frac{1}{C_{\varepsilon_c}} \frac{N}{\log N}\right)$$

$$+ \mathcal{P}_{\varepsilon_c}\left(\left(\frac{\log N}{N}\right)^{5/2} \widetilde{S}_{\lfloor 1/C_{\varepsilon_c} N/\log N \rfloor} > \frac{K}{2}\right).$$

Letting $N \to \infty$, the first term in the second line of this equation vanishes because of (6.27), while for the second term, by (6.23), we have

$$\mathcal{P}_{\varepsilon_c}\left(\left(\frac{\log N}{N}\right)^{5/2} \widetilde{S}_{\lfloor 1/C_{\varepsilon_c} N/\log N \rfloor} > \frac{K}{2}\right) \longrightarrow P\left(\widetilde{L}_1 > \frac{K(C_{\varepsilon_c})^{5/2}}{2}\right).$$

Since $P(\widetilde{L}_1 > t) \to 0$ as $t \to +\infty$, equation (6.31) is proven.

## APPENDIX A: SOME RENEWAL THEORY ESTIMATES

**A.1. Proof of equation (3.18).** We are going to prove equation (3.18), that can be rewritten in terms of the law $\mathcal{P}_\varepsilon$, thanks to (2.18), as

(A.1) $\quad \mathcal{P}_{\varepsilon_c}\left(\delta_N \geq t \frac{N}{\log N} \Big| N+1 \in \chi\right) \leq \frac{c_1}{t} + a_N \quad$ with $a_N \to 0$ as $N \to \infty$,

where we recall that $\delta_N$ has been defined in (3.15). We first need to recall some preliminary relations. We are in the critical case, hence, $q_{\varepsilon_c}(n) = \mathcal{P}_{\varepsilon_c}(\chi_1 = n) \sim C_{\varepsilon_c}/n^2$ by (3.5), because $\mathrm{F}(\varepsilon_c) = 0$. Since $\{\chi_k\}_{k \geq 0}$ is a genuine renewal process, Theorem 8.8.1 of [3] yields

$$\frac{\chi_k}{k \log k} \longrightarrow C_{\varepsilon_c} \qquad \text{as } k \to \infty, \mathcal{P}_{\varepsilon_c}\text{-a.s.}$$

By the definition (3.3) of the variable $\iota_N$, we have $\chi_{\iota_N} \leq N \leq \chi_{\iota_N+1}$, hence,

(A.2) $\qquad \dfrac{\iota_N}{N/\log N} \longrightarrow \dfrac{1}{C_{\varepsilon_c}} \qquad$ as $N \to \infty, \mathcal{P}_{\varepsilon_c}$-a.s.

Introducing the renewal function $u_{\varepsilon_c}(n) := \mathcal{P}_{\varepsilon_c}(n \in \chi) = \sum_{k=0}^{\infty}(q_{\varepsilon_c})^{*k}(n)$, Theorem 8.7.4 of [3] gives

(A.3) $\qquad u_{\varepsilon_c}(n) \sim \dfrac{1}{C_{\varepsilon_c} \log n} \qquad$ as $n \to \infty,$

which implies

(A.4) $\qquad U_{\varepsilon_c}(n) := \sum_{k=0}^{n} u_{\varepsilon_c}(k) \sim \dfrac{n}{C_{\varepsilon_c} \log n} \qquad$ as $n \to \infty.$

We are ready to prove (A.1). We denote by $\xi$ the length of the excursion of $\chi$ embracing the point $N/2$:

(A.5) $$\xi := \min\{\chi \cap [N/2, \infty)\} - \max\{\chi \cap [0, N/2]\}.$$

Then the inclusion bound and the symmetry $n \mapsto N - n$ yield

(A.6) $$\begin{aligned}&\mathcal{P}_{\varepsilon_c}\left(\delta_N \geq t \frac{N}{\log N} \Big| N + 1 \in \chi\right) \\ &\leq \mathcal{P}_{\varepsilon_c}\left(\xi \geq \frac{N}{\log N} \Big| N + 1 \in \chi\right) \\ &\quad + 2\mathcal{P}_{\varepsilon_c}\left(\delta_{\lfloor N/2 \rfloor} \geq t \frac{N}{\log N}, \xi < \frac{N}{\log N} \Big| N + 1 \in \chi\right).\end{aligned}$$

Let us focus on the first term in the r.h.s. of (A.6). We can write

(A.7) $$\mathcal{P}_{\varepsilon_c}\left(\xi \geq \frac{N}{\log N} \Big| N + 1 \in \chi\right) = \sum_{\substack{0 \leq i \leq N/2 \leq j \leq N+1 \\ j - i \geq N/\log N}} \frac{u(i) q(j-i) u(N+1-j)}{u(N+1)},$$

where we have omitted for simplicity the dependence of $q(\cdot)$ and $u(\cdot)$ on $\varepsilon_c$. If we consider the terms in the sum with $i \leq N/4$, then $j - i \geq N/4$ and, therefore, $q(j - i) \leq (const.)/N^2$, hence, recalling (A.3) and (A.4), the contribution of these terms is bounded above by

(A.8) $$\frac{(const.)}{N^2} \frac{U(\lfloor N/2 \rfloor)^2}{u(N+1)} \leq \frac{(const.')}{\log N}.$$

By symmetry, the same bound holds for the contribution of the terms with $j \geq 3N/4$. It remains to consider the terms where both $i > N/4$ and $j < 3N/4$: applying (A.3) to $u(N)$, $u(i)$ and $u(N + 1 - j)$, the contribution of these terms is bounded above by

(A.9) $$\begin{aligned}\frac{(const.)}{(\log N)} &\sum_{\substack{N/4 \leq i \leq N/2 \leq j \leq 3N/4 \\ j - i \geq N/\log N}} \frac{1}{(j-i)^2} \leq \frac{(const.)}{(\log N)} \sum_{l = \lceil N/\log N \rceil}^{\lfloor N/2 \rfloor} (l+1) \cdot \frac{1}{l^2} \\ &\leq (const.') \frac{\log \log N}{\log N}.\end{aligned}$$

We have thus shown that the first term in the r.h.s. of (A.6) vanishes as $N \to \infty$, hence, it can be absorbed in the term $a_N$, appearing in the r.h.s. of (A.1).

Next we consider the second term in the r.h.s. of (A.6). We sum over the location $m$ of $\chi_{\iota_{\lfloor N/2 \rfloor}}$, that is, the last point of $\chi$ before $\lfloor N/2 \rfloor$, and over the location $l$ of $\chi_{\iota_{\lfloor N/2 \rfloor}+1}$, that is, the first point of $\chi$ after $\lfloor N/2 \rfloor$. Recalling (3.15), for $t > 1$ the renewal property yields

$$
\begin{aligned}
&\mathcal{P}_{\varepsilon_c}\left(\delta_{\lfloor N/2 \rfloor} \geq t \frac{N}{\log N}, \xi < \frac{N}{\log N} \Big| N+1 \in \chi\right) \\
&\quad = \sum_{\substack{m \leq \lfloor N/2 \rfloor, l > \lfloor N/2 \rfloor \\ l-m < N/\log N}} \mathcal{P}_{\varepsilon_c}\left(\delta_m \geq t \frac{N}{\log N}, m \in \chi\right) \\
&\qquad \times q(l-m) \cdot \frac{u(N+1-l)}{u(N+1)}.
\end{aligned}
$$
(A.10)

In the range of summation, by (A.3), the ratio $u(N+1-l)/u(N+1)$ is bounded above by some positive constant $A$, hence, the r.h.s. is bounded above by

$$
A \sum_{\substack{m \leq \lfloor N/2 \rfloor, l > \lfloor N/2 \rfloor \\ l-m < N/\log N}} \mathcal{P}_{\varepsilon_c}\left(\delta_m \geq t \frac{N}{\log N}, m \in \chi\right) q(l-m)
$$

$$
\leq A \mathcal{P}_{\varepsilon_c}\left(\delta_{\lfloor N/2 \rfloor} \geq t \frac{N}{\log N}\right).
$$

We are finally reduced to estimating the last term. By (A.2), we can write as $N \to \infty$

$$
\mathcal{P}_{\varepsilon_c}\left(\delta_{\lfloor N/2 \rfloor} \geq t \frac{N}{\log N}\right) = \mathcal{P}_{\varepsilon_c}\left(\delta_{\lfloor N/2 \rfloor} \geq t \frac{N}{\log N}, \iota_{\lfloor N/2 \rfloor} \leq \frac{2}{C_{\varepsilon_c}} \frac{N/2}{\log N}\right) + o(1)
$$

and, by (3.15), the first term in the r.h.s. is bounded above by

$$
\mathcal{P}_{\varepsilon_c}\left(\max\left\{\chi_i - \chi_{i-1} : i \leq \frac{N}{C_{\varepsilon_c} \log N}\right\} \geq t \frac{N}{\log N}\right).
$$

This probability is easily estimated. In fact, the variables $\{\chi_i - \chi_{i-1}\}_{i \in \mathbb{N}}$ under $\mathcal{P}_{\varepsilon_c}$ are independent and identically distributed, hence, for $x > 0$ and $M \in \mathbb{N}$ we have

$$
\mathcal{P}_{\varepsilon_c}(\max\{\chi_i - \chi_{i-1} : i \leq M\} < x) = \mathcal{P}_{\varepsilon_c}(\chi_1 < x)^M \geq \left(1 - \frac{B}{x}\right)^M,
$$

where $B$ is a suitable positive constant. Since $(1-t) \geq e^{-2t}$ for $t \in [0, \frac{1}{2}]$, it follows that for $N$ sufficiently large we have

$$
\mathcal{P}_{\varepsilon_c}\left(\max\left\{\chi_i - \chi_{i-1} : i \leq \frac{N}{C_{\varepsilon_c} \log N}\right\} \geq t \frac{N}{\log N}\right) \leq 1 - \exp\left(-\frac{2B}{C_{\varepsilon_c} t}\right) \leq \frac{2B}{C_{\varepsilon_c} t}
$$

and the proof of (A.1) is completed.

**A.2. Proof of equation (5.1).** In this section we prove (5.1), which we can rewrite as

$$(A.11) \qquad \lim_{t \to 0^+} \liminf_{N \to \infty} \mathcal{P}_{\varepsilon_c}\left(\delta_N \geq t \frac{N}{\log N} \Big| N+1 \in \chi\right) = 1.$$

We start observing that the inclusion bound yields

$$\mathcal{P}_{\varepsilon_c}\left(\delta_N \geq t \frac{N}{\log N} \Big| N+1 \in \chi\right) \geq \mathcal{P}_{\varepsilon_c}\left(\delta_{\lfloor N/2 \rfloor} \geq t \frac{N}{\log N}, \xi < \frac{N}{\log N} \Big| N+1 \in \chi\right),$$

where we recall that the variable $\xi$ has been defined in (A.5). We decompose the r.h.s. according to (A.10) and we observe that the fraction $u(N+1-l)/u(N+1)$ converges to 1 as $N \to \infty$ uniformly in the range of summation, by (A.3). Therefore, we can write

$$\mathcal{P}_{\varepsilon_c}\left(\delta_N \geq t \frac{N}{\log N}, \xi < \frac{N}{\log N} \Big| N+1 \in \chi\right)$$

$$\geq (1+o(1)) \sum_{\substack{m \leq \lfloor N/2 \rfloor, l > \lfloor N/2 \rfloor \\ l-m < N/\log N}} \mathcal{P}_{\varepsilon_c}\left(\delta_m \geq t \frac{N}{\log N}, m \in \chi\right) q_{\varepsilon_c}(l-m)$$

$$= (1+o(1)) \mathcal{P}_{\varepsilon_c}\left(\delta_{\lfloor N/2 \rfloor} \geq t \frac{N}{\log N}, \xi < \frac{N}{\log N}\right) \qquad (N \to \infty).$$

Recalling that $q_{\varepsilon_c}(n) \sim C_{\varepsilon_c}/n^2$ and $u_{\varepsilon_c}(n) \sim 1/(C_{\varepsilon_c} \log n)$ as $n \to \infty$, by (3.5) and (A.3), we obtain

$$(A.12) \quad \mathcal{P}_{\varepsilon_c}\left(\xi \geq \frac{N}{\log N}\right) = \sum_{\substack{m \leq \lfloor N/2 \rfloor, l > \lfloor N/2 \rfloor \\ l-m \geq N/\log N}} u_{\varepsilon_c}(m) q_{\varepsilon_c}(l-m) \leq \frac{(const.)}{\log N}.$$

Putting together the preceding relations, we have

$$\mathcal{P}_{\varepsilon_c}\left(\delta_N \geq t \frac{N}{\log N} \Big| N+1 \in \chi\right) \geq \mathcal{P}_{\varepsilon_c}\left(\delta_{\lfloor N/2 \rfloor} \geq t \frac{N}{\log N}\right) + o(1) \qquad (N \to \infty)$$

and we are left with estimating the r.h.s. of this relation. The inclusion bound, the definition (3.15) of $\delta_N$ and equation (A.2) yield

$$\mathcal{P}_{\varepsilon_c}\left(\delta_{\lfloor N/2 \rfloor} \geq t \frac{N}{\log N}\right)$$

$$\geq \mathcal{P}_{\varepsilon_c}\left(\delta_{\lfloor N/2 \rfloor} \geq t \frac{N}{\log N}, \iota_{\lfloor N/2 \rfloor} \geq \frac{1}{2C_{\varepsilon_c}} \frac{N/2}{\log N}\right)$$

$$\geq \mathcal{P}_{\varepsilon_c}\left(\max\left\{\chi_i - \chi_{i-1} : i \leq \frac{1}{4C_{\varepsilon_c}} \frac{N}{\log N}\right\} \geq t \frac{N}{\log N}, \iota_{\lfloor N/2 \rfloor} \geq \frac{1}{2C_{\varepsilon_c}} \frac{N/2}{\log N}\right)$$

$$= \mathcal{P}_{\varepsilon_c}\left(\max\left\{\chi_i - \chi_{i-1} : i \leq \frac{1}{4C_{\varepsilon_c}} \frac{N}{\log N}\right\} \geq t \frac{N}{\log N}\right) - o(1) \qquad (N \to \infty).$$

The variables $\{\chi_i - \chi_{i-1}\}_{i \in \mathbb{N}}$ under $\mathcal{P}_{\varepsilon_c}$ are independent and identically distributed, hence, for $x > 0$ and $M \in \mathbb{N}$ we have

$$\mathcal{P}_{\varepsilon_c}(\max\{\chi_i - \chi_{i-1} : i \leq M\} < x) = \mathcal{P}_{\varepsilon_c}(\chi_1 < x)^M$$
$$\leq \left(1 - \frac{D}{x}\right)^M \leq e^{-D/xM}$$

for some positive constant $D$. Therefore,

$$\mathcal{P}_{\varepsilon_c}\left(\delta_{\lfloor N/2 \rfloor} \geq t\frac{N}{\log N}\right) \geq 1 - \exp\left(-\frac{D}{4C_{\varepsilon_c}t}\right) + o(1) \qquad (N \to \infty)$$

and the proof of relation (A.11) is complete.

**A.3. Proof of equation (3.19).** We are going to prove equation (3.19), which can be rewritten using (2.18) as

(A.13)    $\mathcal{P}_\varepsilon(\delta_N \geq c_2 \log N | N + 1 \in \chi) \longrightarrow 0 \qquad$ as $N \to \infty$.

Since we assume that $\varepsilon > \varepsilon_c$, we are in the localized regime and the step law $q_\varepsilon(n) = \mathcal{P}_\varepsilon(\chi_1 = n)$ has exponential tails; see (3.6). The renewal theorem then yields

(A.14)    $\mathcal{P}_\varepsilon(N \in \chi) \longrightarrow \frac{1}{\mathcal{E}_\varepsilon(\chi_1)} \in (0, \infty) \qquad$ as $N \to \infty$

and the weak law of large numbers gives

$$\mathcal{P}_\varepsilon\left(\iota_N \geq \frac{2}{\mathcal{E}_\varepsilon(\chi_1)}N\right) \longrightarrow 0 \qquad \text{as } N \to \infty.$$

These relations yield

$$\mathcal{P}_\varepsilon(\delta_N \geq c_2 \log N | N + 1 \in \chi)$$
$$\leq (const.)\mathcal{P}_\varepsilon(\delta_N \geq c_2 \log N)$$
$$= (const.)\mathcal{P}_\varepsilon\left(\delta_N \geq c_2 \log N, \iota_N \leq \frac{2N}{\mathcal{E}_\varepsilon(\chi_1)}\right) + o(1) \qquad \text{as } N \to \infty.$$

The definition (3.15) of $\delta_N$ and the inclusion bound give

$$\mathcal{P}_\varepsilon\left(\delta_N \geq c_2 \log N, \iota_N \leq \frac{2N}{\mathcal{E}_\varepsilon(\chi_1)}\right)$$
$$\leq \mathcal{P}_\varepsilon\left(\max\left\{\chi_i - \chi_{i-1} : i \leq \frac{2N}{\mathcal{E}_\varepsilon(\chi_1)}\right\} \geq c_2 \log N\right).$$

Since the variables $\{\chi_i - \chi_{i-1}\}_{i \in \mathbb{N}}$ under $\mathcal{P}_\varepsilon$ are independent and identically distributed, for $x > 0$ and $M \in \mathbb{N}$ we have

$$\mathcal{P}_\varepsilon(\max\{\chi_i - \chi_{i-1} : i \leq M\} < x) = \mathcal{P}_\varepsilon(\chi_1 < x)^M \geq (1 - Be^{-\mathrm{G}(\varepsilon)x})^M$$

for a suitable positive constant $B$. Since $(1-t) \geq e^{-2t}$ for $t \in [0, \frac{1}{2}]$, it follows that for $N$ sufficiently large we have

$$\mathcal{P}_\varepsilon\left(\delta_N \geq c_2 \log N, \iota_N \leq \frac{2N}{\mathcal{E}_\varepsilon(\chi_1)}\right) \leq 1 - \exp\left(-2B \frac{2N}{\mathcal{E}_\varepsilon(\chi_1)} \frac{1}{N^{\mathrm{G}(\varepsilon) \cdot c_2}}\right).$$

If we choose $c_2 > 1/\mathrm{G}(\varepsilon)$, the r.h.s. vanishes as $N \to \infty$ and equation (A.13) is proven.

## APPENDIX B: SOME TECHNICAL PROOFS

**B.1. Proof of Proposition 6.3.** Take any subsequence $\{\boldsymbol{\nu}_{N_n}\}_{n \in \mathbb{N}}$ that converges in distribution toward some random signed measure $\boldsymbol{\nu}$. We are going to show that the finite dimensional distributions of $\boldsymbol{\nu}$ are necessarily given by the laws $\lambda^{(k)}_{a_1, \ldots, a_k}$ that appear in (6.3). Since the finite-dimensional distributions determine laws on $\mathcal{M}([0,1])$, this means that every convergent subsequence of $\{\boldsymbol{\nu}_N\}_{N \in \mathbb{N}}$ must have the same limit. Then Lemma 6.2 and a standard sub-subsequence argument yield the convergence of the whole sequence $\{\boldsymbol{\nu}_N\}_{N \in \mathbb{N}}$, and the proof is complete.

Therefore, we assume that $\{\boldsymbol{\nu}_{N_n}\}_{n \in \mathbb{N}}$ converges in distribution toward $\boldsymbol{\nu}$. We introduce the function $f_t^{(\varepsilon)} : [0,1] \to \mathbb{R}$ defined by

$$f_t^{(\varepsilon)}(x) := \begin{cases} 1, & x \in [0, t], \\ -\dfrac{x}{\varepsilon} + 1 + \dfrac{t}{\varepsilon}, & x \in [t, t + \varepsilon], \\ 0, & x \in [t + \varepsilon, 1], \end{cases}$$

which may be viewed as a continuous approximation of $\mathbf{1}_{[0,t]}$. Then we define the map $F_t^{(\varepsilon)} : \mathcal{M}([0,1]) \to \mathbb{R}$ by $F_t^{(\varepsilon)}(\nu) := \int f_t^{(\varepsilon)} \, d\nu$. Notice that $|F_t^{(\varepsilon)}(\nu) - F_t(\nu)| \leq |\nu|([t, t+\varepsilon])$. Now let $W : \mathbb{R}^k \to \mathbb{R}$ be a bounded and Lipschitz function such that

(B.1)
$$|W(x_1, \ldots, x_k) - W(y_1, \ldots, y_k)| \leq \sum_{i=1}^k g(x_i - y_i)$$

where $g(x) := |x| \wedge 1$.

Therefore, we can write

(B.2)
$$|E[W(F_{a_1}^{(\varepsilon)}(\boldsymbol{\nu}_{N_n}), \ldots, F_{a_k}^{(\varepsilon)}(\boldsymbol{\nu}_{N_n}))] - E[W(F_{a_1}(\boldsymbol{\nu}_{N_n}), \ldots, F_{a_k}(\boldsymbol{\nu}_{N_n}))]|$$
$$\leq \sum_{i=1}^k E[g(|\boldsymbol{\nu}_{N_n}|([a_i, a_i + \varepsilon]))].$$

Let us take the $n \to \infty$ limit. Since $W(\cdot)$ and $F_t^{(\varepsilon)}(\cdot)$ are continuous,

$$E[W(F_{a_1}^{(\varepsilon)}(\boldsymbol{\nu}_{N_n}), \ldots, F_{a_k}^{(\varepsilon)}(\boldsymbol{\nu}_{N_n}))] \longrightarrow E[W(F_{a_1}^{(\varepsilon)}(\boldsymbol{\nu}), \ldots, F_{a_k}^{(\varepsilon)}(\boldsymbol{\nu}))],$$

and also $E[W(F_{a_1}(\boldsymbol{\nu}_{N_n}), \ldots, F_{a_k}(\boldsymbol{\nu}_{N_n}))] \longrightarrow \int W \, \mathrm{d}\lambda_{a_1,\ldots,a_k}^{(k)}$ by (6.3). Then we take the limit $\varepsilon \to 0$: the r.h.s. of (B.2) vanishes by (6.4) and by dominated convergence, we have

$$E[W(F_{a_1}(\boldsymbol{\nu}), \ldots, F_{a_k}(\boldsymbol{\nu}))] = \int W \, \mathrm{d}\lambda_{a_1,\ldots,a_k}^{(k)}.$$

Since $W(\cdot)$ is an arbitrary function satisfying (B.1), this shows that the finite dimensional distributions of $\boldsymbol{\nu}$ are indeed $\lambda_{a_1,\ldots,a_k}^{(k)}$, and the proof is complete.

**B.2. Proof of Lemma 6.2.** Let us denote by $\nu_N := \boldsymbol{\nu}_N \circ P^{-1}$ the law of the random signed measure $\boldsymbol{\nu}_N$, so that $\nu_N$ is a probability measure on $\mathcal{M}([0,1])$. For every fixed $K \in \mathbb{N}$, the restriction of $\nu_N$ on the subspace $\mathcal{M}_K([0,1])$ is a sub-probability measure on a Polish space (cf. Lemma 6.1), hence, one can apply the standard Prohorov theorem. So we can extract a subsequence $\{\nu_{N'}\}$ that converges weakly toward a sub-probability law $\lambda^{(1)}$ on $\mathcal{M}_1([0,1])$; then from $\{\nu_{N'}\}$ we extract a sub-subsequence $\{\nu_{N''}\}$ that converges weakly toward a sub-probability law $\lambda^{(2)}$ on $\mathcal{M}_2([0,1])$, and so on. With a standard diagonal argument, we obtain a subsequence $\{\nu_{N_k}\}_k$ that converges weakly on $\mathcal{M}_K([0,1])$ toward $\lambda^{(K)}$, for every $K \in \mathbb{N}$. However, recalling (6.1), it is clear that the laws $\lambda^{(K)}$ are the restriction on $\mathcal{M}_K([0,1])$ of a single law $\lambda$ on $\mathcal{M}([0,1])$ and, moreover, $\lambda(\mathcal{M}([0,1])) = 1$ because the sequence $\{\boldsymbol{\nu}_N\}_N$ is tight; compare (6.2). Then it is easy to check that the subsequence $\{\nu_{N_k}\}_k$ converges weakly on $\mathcal{M}([0,1])$ toward $\lambda$: in fact, given a continuous and bounded functional $G: \mathcal{M}([0,1]) \to \mathbb{R}$, we can write

$$\left| \int G \, \mathrm{d}\nu_{N_k} - \int G \, \mathrm{d}\lambda \right|$$
$$\leq \left| \int G \mathbf{1}_{\mathcal{M}_K([0,1])} \, \mathrm{d}\nu_{N_k} - \int G \mathbf{1}_{\mathcal{M}_K([0,1])} \, \mathrm{d}\lambda \right|$$
$$+ \|G\|_\infty \cdot (\nu_{N_k}(\mathcal{M}_K([0,1])^\complement) + \lambda(\mathcal{M}_K([0,1])^\complement)).$$

The first term in the r.h.s. vanishes as $k \to \infty$, because, by construction, $\nu_{N_k}$ converges weakly to $\lambda = \lambda^{(K)}$ on $\mathcal{M}_K([0,1])$, and the second term vanishes as $K \to \infty$ because of the tightness of $\{\boldsymbol{\nu}_N\}_N$; compare (6.2). This completes the proof.

**B.3. Computing $\Phi(t)$.** We recall that $\{B_t\}_{t \in [0,1]}$ is a standard Brownian motion on $\mathbb{R}$ and $I_t := \int_0^t B_s \, \mathrm{d}s$. We also set $G_t := \int_0^t I_s \, \mathrm{d}s$. The function $\Phi(t)$ was introduced in (6.14): recalling the definition (1.11) of the conditioned

44	F. CARAVENNA AND J.-D. DEUSCHELprocess $\widehat{I}_t$, we can re-express it as

$$\Phi(t) = P(G_1 > t | B_1 = 0, I_1 = 0).$$

Since the vector $(G_1, I_1, B_1)$ has a centered Gaussian distribution, the law of $G_1$ under $P(\cdot | B_1 = 0, I_1 = 0)$ is centered Gaussian too and, hence, it suffices to identify its variance to determine $\Phi(t)$. The covariance matrix $A$ of $(G_1, I_1, B_1)$ is easily computed:

$$A := \begin{pmatrix} E(G_1^2) & E(G_1 I_1) & E(G_1 B_1) \\ E(G_1 I_1) & E(I_1^2) & E(I_1 B_1) \\ E(G_1 B_1) & E(I_1 B_1) & E(B_1^2) \end{pmatrix} = \begin{pmatrix} \frac{1}{20} & \frac{1}{8} & \frac{1}{6} \\ \frac{1}{8} & \frac{1}{3} & \frac{1}{2} \\ \frac{1}{6} & \frac{1}{2} & 1 \end{pmatrix}.$$

The variance of $G_1$ conditionally on $\{I_1 = 0, B_1 = 0\}$ is then given by $1/(A^{-1})_{1,1} = \frac{1}{720}$. Therefore,

$$\Phi(t) = \int_t^\infty \frac{e^{-360s^2}}{\sqrt{2\pi/720}} \, ds = \frac{6\sqrt{10}}{\sqrt{\pi}} \int_t^\infty e^{-360s^2} \, ds.$$

**Acknowledgments.** We are very grateful to Yvan Velenik and Ostap Hryniv for suggesting the present problem to us.## REFERENCES

[1] ASMUSSEN, S. (2003). *Applied Probability and Queues*, 2nd ed. *Applications of Mathematics: Stochastic Modelling and Applied Probability (New York)* **51**. Springer, New York. MR1978607
[2] BERTOIN, J. (1996). *Lévy Processes. Cambridge Tracts in Mathematics* **121**. Cambridge Univ. Press, Cambridge. MR1406564
[3] BINGHAM, N. H., GOLDIE, C. M. and TEUGELS, J. L. (1989). *Regular Variation. Encyclopedia of Mathematics and Its Applications* **27**. Cambridge Univ. Press, Cambridge. MR1015093
[4] BRASCAMP, H. J. and LIEB, E. H. (1976). On extensions of the Brunn–Minkowski and Prékopa–Leindler theorems, including inequalities for log concave functions, and with an application to the diffusion equation. *J. Functional Analysis* **22** 366–389. MR0450480
[5] CAFFARELLI, L. A. (2000). Monotonicity properties of optimal transportation and the FKG and related inequalities. *Comm. Math. Phys.* **214** 547–563. MR1800860
[6] CARAVENNA, F. and DEUSCHEL, J.-D. (2008). Pinning and wetting transition for $(1 + 1)$-dimensional fields with Laplacian interaction. *Ann. Probab.* **36** 2388–2433.
[7] COHN, D. L. (1980). *Measure Theory*. Birkhäuser, Boston, MA. MR578344
[8] FELLER, W. (1971). *An Introduction to Probability Theory and Its Applications* **II**, 2nd ed. Wiley, New York. MR0270403
[9] HÖGNÄS, G. (1977). Characterization of weak convergence of signed measures on $[0, 1]$. *Math. Scand.* **41** 175–184. MR0482909
[10] JOHNSON, J. (1985). An elementary characterization of weak convergence of measures. *Amer. Math. Monthly* **92** 136–137. MR777561
[11] KINGMAN, J. F. C. (1993). *Poisson Processes. Oxford Studies in Probability* **3**. Oxford Univ. Press, New York. MR1207584

Dipartimento di Matematica Pura e Applicata  
Università degli Studi di Padova  
via Trieste 63  
35121 Padova  
Italy  
E-mail: francesco.caravenna@math.unipd.it

Institut für Mathematik  
TU Berlin  
Strasse des 17, Juni 136  
10623 Berlin  
Germany  
E-mail: deuschel@math.tu-berlin.de